\documentclass[a4paper]{article}

\usepackage{amsmath}

\usepackage{amsthm}
\usepackage[dvips]{graphicx}

\usepackage{epsfig}

\usepackage{amssymb}
\usepackage{latexsym}
\newtheorem{Definition}{Definition}
\newtheorem{Theorem}{Theorem}
\newtheorem{Proposition}{Proposition}
\newtheorem{Lemma}{Lemma}
\newtheorem{Corollary}{Corollary}
\newtheorem{Remark}{Remark}
\newtheorem{Claim}{Claim}

\def\demo{ {\bf Proof.} }

\def\SS{{\mathcal S}}

\def\II{{\mathcal I}}

\def\OO{{\mathcal O}}

\def\KK{{\mathcal K}}

\def\A{{\mathcal A}}

\def\UU{{\mathcal U}}
\def\KK{{\mathcal K}}

\def\Z{{\mathbb Z}}
\def\Q{{\mathbb Q}}

\def\S{{\mathbb S}}

\begin{document}

\title{On cyclic branched coverings of prime knots}

\author{Michel Boileau and Luisa Paoluzzi} 

\date{\today}

\maketitle

\bigskip 
\bigskip

\begin{abstract}
We prove that a prime knot $K$ is not determined by its $p$-fold cyclic 
branched cover for at most two odd primes $p$. Moreover, we show that for a 
given odd prime $p$, the $p$-fold cyclic branched cover of a prime knot $K$ is
the $p$-fold cyclic branched cover of at most one more knot $K'$ non equivalent
to $K$. To prove the main theorem, a result concerning the symmetries of knots
is also obtained. This latter result can be interpreted as a characterisation
of the trivial knot.

\vskip 2mm


\vskip 2mm


\end{abstract}

\section{Introduction}

Two knots $K$ and $K'$ are \emph{equivalent} if there is a homeomorphism of 
$\S^{3}$ sending $K$ to $K'$. Given a knot $K \subset \S^{3}$ and an integer 
$p \geq 2$ one can construct the (total space of the) $p$-fold cyclic cover 
$M_{p}(K)$ of $\S^{3}$ branched along $K$: it is a fundamental object in knot 
theory. There are non-prime knots all of whose cyclic branched covers are 
homeomorphic. This is no longer true for prime knots: S. Kojima \cite{Ko} 
proved that for each prime knot $K \subset \S^{3}$ there is an integer $n_{K} 
\geq 2$ such that two prime knots $K$ and $K'$ are equivalent if their $p$-fold 
cyclic branched covers are homeomorphic for some $p > \max (n_{K}, n_{K'})$.

\bigskip

There are many examples  of prime knots in $\S^{3}$ which are not equivalent 
but share homeomorphic $p$-fold cyclic branched covers due to C. Giller 
\cite{Gi}, C. Livingston \cite{Li}, Y. Nakanishi \cite{Na}, M. Sakuma 
\cite{Sa1}. Moreover there is no universal bound for $n_{K}$.

The main goal of this article is to study the relationship between prime knots 
and their cyclic branched covers when the number of sheets is an odd prime 
number. 

\smallskip

\begin{Definition} 
Let $K \subset \S^3$ be a prime knot. A knot $K' \subset \S^3$ which is not 
equivalent to $K$ and which has the same $p$-fold cyclic branched cover as $K$ 
is called a \emph{$p$-twin} of $K$.
\end{Definition}

\medskip

There are examples of prime knots, even hyperbolic knots (e.g. Montesinos 
knots) with an arbitrarily large number of non-equivalent $2$-twins. In 
contrast, for an odd prime number $p$, the number of $p$-twins is very 
restricted, according to our main result:

\smallskip

\begin{Theorem}\label{thm:twins} 
Let $K\subset \S^3$ be a prime knot. Then:

\item{(i)} There are at most two odd prime numbers $p$ for which $K$ admits a 
$p$-twin.

\item{(ii)} For a given odd prime number $p$, $K$ admits at most one $p$-twin.

\item{(iii)} Suppose that a prime knot $K$ admits the same knot $K'$ as a 
$p$-twin and a $q$-twin for two distinct odd prime numbers $p$ and $q$. Then 
$K$ has two commuting rotational symmetries of order $p$ and $q$ with trivial 
quotients.
\end{Theorem}

\medskip

A \emph{rotational symmetry of order $p$} of a knot $K \subset \S^3$ is an 
orientation preserving periodic diffeomorphism $\psi$ of the pair $(\S^3, K)$ 
with period $p$ and non-empty fixed-point set disjoint from $K$. We say that 
the rotational symmetry $\psi$ has \emph{trivial quotient} if $K/\psi$ is the 
trivial knot.

For hyperbolic knots Theorem \ref{thm:twins} is in fact a consequence of B. 
Zimmermann's result in \cite{Zim1} whose proof uses the orbifold theorem and 
the Sylow theory for finite groups. 

The result in Theorem \ref{thm:twins} is sharp: for any pair of coprime 
integers $p> q >2$ B. Zimmermann has constructed examples of prime hyperbolic 
knots with the same $p$-fold and $q$-fold branched coverings \cite{Zim2}. 

The second named author \cite{Pao2} has proved that a hyperbolic knot is 
determined by three cyclic branched covers of pairwise distinct orders. The 
following, straightforward corollary of Theorem \ref{thm:twins}, shows that a 
stronger conclusion holds for arbitrary prime knots when we focus on branched 
coverings with odd prime orders.

\smallskip

\begin{Corollary}\label{cor: three covers} A prime knot is determined by three 
cyclic branched covers of pairwise distinct odd prime orders. More 
specifically, for every knot $K$ there is at least one integer $p_K \in 
\{3, 5, 7 \}$ such that $K$ is determined by its $p_K$-cyclic branched cover. 
\end{Corollary}

\medskip

Another straightforwards consequence of Theorem \ref{thm:twins} is:

\smallskip

\begin{Corollary}\label{cor:composite}
Let $K=K_1\sharp...\sharp K_t$ and $K'=K'_1\sharp...\sharp K'_t$ be two 
composite knots with the same cyclic branched covers of orders $p_j$, 
$j=1,2,3$, for three fixed, pairwise distinct, odd prime numbers. Then, after
a reordering, the (non oriented) knots $K_i$ and $K'_i$ are equivalent for all
$i=1,...,t$.
\end{Corollary}

\medskip

Part (ii) of Theorem \ref{thm:twins} states that for a given odd prime number 
$p$ a closed, orientable $3$-manifold can be the $p$-fold cyclic branched cover 
of at most two non-equivalent knots in $\S^3$. In \cite{BPZ} it has been shown 
that an integer homology sphere which is a $n$-fold cyclic branched cover of  
$\S^3$ for four distinct odd prime numbers $n$ is in fact $\S^3$. By putting 
together these two results we get the following corollary:

\begin{Corollary}\label{cor:homologysphere} 
Let $M$ be an irreducible integer homology $3$-sphere. Then: there are at most 
three distinct knots in $\S^3$ having $M$ as cyclic branched cover of odd prime 
order.
\end{Corollary}

Our main task will be to prove Theorem \ref{thm:twins} for a satellite knot: 
that is a knot whose exterior $\S^3\setminus\UU(K)$ has a non trivial 
Jaco-Shalen-Johannson decomposition \cite{JS}, \cite{Jo} (in the sequel we use 
$JSJ$-decomposition for short). Otherwise the knot is called simple: in this 
case, due to Thurston's hyperbolization theorem \cite{Th2}, its exterior is 
either hyperbolic, and the proof follows already from the works in \cite{Pao2} 
and \cite{Zim1}, or it is a torus knot and a simple combinatorial argument 
applies.
 
The proof of Theorem \ref{thm:twins} for satellite knots relies on the study 
of the \emph{partial symmetries} of the exterior $E(K)$ of $K$ induced by the 
covering transformations associated to the twins of $K$ and on the localization 
of their axes of fixed points in the components of the $JSJ$-decomposition of 
$E(K)$. In particular the proof uses the following result about rotational 
symmetries of prime knots which is of interest in its own right.

\smallskip

\begin{Theorem}\label{thm:three rotations} 
Let $K$ be a knot in $\S^3$ admitting three rotational symmetries  with trivial 
quotients and whose orders are three pairwise distinct numbers $>2$. Then $K$ 
is the trivial knot.
\end{Theorem}

\medskip

Since the trivial knot admits a rotational symmetry with trivial quotient of
order $p$ for each integer $p \ge 2$, the above Theorem \ref{thm:three
rotations} can be interpreted as a characterisation of the trivial knot, i.e. a
knot is trivial if and only if it admits three rotational symmetries of
pairwise distinct orders $>2$ and trivial quotients.

\section{Rotational symmetries of knots}

A \emph{rotational symmetry} of order $p$ of a knot $K \subset \S^3$ is an 
orientation preserving, periodic diffeomorphism $\psi$ of the pair $(\S^3, K)$ 
of order $p$ and non-empty fixed-point set disjoint from $K$. We say that the 
rotational symmetry $\psi$ has \emph{trivial quotient} if $K/\psi$ is the 
trivial knot.

\smallskip

\begin{Remark}\label{rem:lift} 
Let $K$ be a knot and let $\psi$ be a rotational symmetry of $K$ of order $p$. 
The symmetry $\psi$ lifts to a periodic diffeomorphism $\tilde\psi$ of the 
$p$-fold branched cover $M_p(K)$ with order $p$ and non-empty fixed-point set, 
which commutes with the covering transformation $h$ of $K$ acting on $M_p(K)$. 
Then the symmetry $\psi$ has trivial quotient if and only if 
$(M,Fix(\tilde\psi))/<\tilde\psi> \cong (\S^3,K')$. Moreover in this case $K$ 
and $K'$ have a common quotient link with two trivial components (see 
\cite{Zim1}).

In particular a symmetry of a knot $K$ induced by the covering transformation 
associated to a $p$-twin $K'$ of $K$ is a $p$-rotational symmetry with trivial 
quotient. This follows from the fact that the two commuting deck 
transformations associated to the two twins induce on $M_p(K)$ a 
$\Z/p\Z \oplus \Z/p\Z$-cover of $\S^3$ branched over a link with two unknotted 
components.
\end{Remark}

\medskip

The main result of this section is the following theorem whose assertion (i) is 
Theorem \ref{thm:three rotations}:

\smallskip

\begin{Theorem}\label{thm:rotations} 
Let $K$ be a knot in $\S^3$.

\item{(i)} Assume that $K$ admits three rotational symmetries  with trivial 
quotients and whose orders are three pairwise distinct numbers $>2$. Then $K$ 
is the trivial knot.

\item{(ii)} Assume that $K$ admits two rotational symmetries $\psi$ and 
$\varphi$ with trivial quotients and of distinct orders $>2$. Then the 
fixed-point sets $Fix(\psi)$ and $Fix(\varphi)$ sit in the $JSJ$-component of 
$E(K)$ which contains $\partial E(K)$.

\end{Theorem}

\medskip

We prove first a weaker version of Theorem \ref{thm:three rotations} that we 
shall use in the remaining of this section (see also \cite[Scholium]{Pao2}).

\smallskip

\begin{Proposition}\label{prop:commuting rotations} 
Let $K$ be a knot in $\S^3$ admitting three commuting rotational symmetries of
orders $p>q>r\ge2$. If the symmetries of order $q$ and $r$ have trivial
quotients, then $K$ is the trivial knot.
\end{Proposition}

\demo
Denote by $\varphi$, $\psi$ and $\rho$ the three symmetries. If two of them 
-say $\varphi$, $\psi$- have the same axis, then by hypothesis the one with 
smaller order -say $\psi$- must have trivial quotient, i.e. $K/\psi$ is the 
trivial knot. Since the three symmetries commute, $\varphi$ induces a 
rotational symmetry of $K/\psi$ which is non trivial for the order of 
$\varphi$ is larger than that of $\psi$. The axis $\A$ of this induced symmetry 
is the image of $Fix(\psi)$ in the quotient by the action of $\psi$. In 
particular $K/\psi$ and $\A$ form a Hopf link and $K$ is the trivial knot: this 
follows from the equivariant Dehn lemma, see \cite{Hil}. We can thus assume 
that the axes are pairwise disjoint. Note that even if $r=2$, since the 
symmetries commute, the symmetry of order $2$ cannot act as a strong inversion 
on the axes of the other two symmetries. In this case we would have that the 
axis of $\rho$, which is a trivial knot, admits two commuting rotational 
symmetries, $\varphi$ and $\psi$, with distinct axes, which is impossible: this 
follows, for instance, from the fact (see \cite[Thm 5.2]{EL}) that one can find 
a fibration of the complement of the trivial knot which is equivariant with 
respect to the two symmetries.
\qed

\bigskip

The proof of Theorem \ref{thm:rotations} is based on a series of Lemmata. 

The first result concerns the structure of the $JSJ$-decomposition of the 
$p$-fold cyclic branched cover $M$ of a prime knot $K \subset \S^3$. Let $h$ be 
the covering transformation, then the quotient space $M/<h>$ has a natural 
orbifold structure, denoted by $\OO_p(K)$, with underlying space $\S^3$ and  
singular locus $K$ with local group a cyclic group of order $p$ (cf.
\cite[Chap. 2]{BMP}). According to Bonahon-Siebenmann \cite{BS} and the 
orbifold theorem \cite{BoP}, \cite{CHK}, such an orbifold admits a 
characteristic collection of toric $2$-suborbifolds, which split $\OO_p(K)$ 
into geometric suborbifolds. Moreover this characteristic collection of toric 
$2$-suborbifolds lifts to the $JSJ$-collection of tori for $M$. It follows that 
for $p > 2$ the Bonahon-Siebenmann characteristic collection of toric 
$2$-suborbifolds coincides with the $JSJ$-collection of tori for the exterior 
$E(K) = \S^3\setminus\UU(K)$ of $K$.

\smallskip

\begin{Lemma}\label{lem:JSJ}
Let $p >2$ be an integer and let $M$ be the $p$-fold cyclic branched cover of a 
prime knot $K$ in the $3$-sphere. Then:

\item{(a)} The dual graph associated to the $JSJ$-decomposition of $M$ is a 
tree.

\item{(b)} The fixed-point set of the group of deck transformations is entirely 
contained in one geometric piece of the decomposition. 
\end{Lemma}

\demo

\noindent{\bf (a)} Note, first of all, that $M$ is irreducible since $K$ is 
prime. Hence the Bonahon-Siebenmann decomposition of the orbifold $\OO_p(K)$ 
lifts to the $JSJ$-collection for $M$ since $p>2$. Moreover, the graph dual to 
the Bonahon-Siebenmann decomposition of the orbifold $\OO_p(K)$, which lifts to 
the $JSJ$-decomposition for $M$, is a tree. Cutting along a torus of former 
decomposition and considering the component $C$ which does not contain 
$K$ one gets the complement of a knot in $\S^3$. The lemma follows now from the 
fact that each connected component of a cyclic branched cover of $C$ has a
unique boundary component.

\noindent{\bf (b)} Note that the group of deck transformations preserves the 
$JSJ$-collection of tori. If $p>2$, the fixed-point set of this group does not 
meet any torus of the $JSJ$-decomposition, because each $JSJ$-torus is
separating and the fixed point set is connected. Since the fixed point set is 
connected, it is entirely contained in one geometric piece of the 
$JSJ$-decomposition.
\qed

\bigskip

\begin{Remark}
Note that the conclusion of the first part of the lemma holds also
for covers of order $2$. For covers of prime order this
property follows also from the fact that $M_p(K)$ is a $\Z/p\Z$-homology sphere (see
\cite{Go}). 
\end{Remark}

\medskip

\begin{Lemma}\label{lem:prime} 
If a knot $K \subset \S^3$ has a rotational symmetry with trivial quotient, 
then $K$ is prime.
\end{Lemma}

\demo 
M. Sakuma \cite[Thm 4]{Sa2} showed that the only possible rotational symmetries 
of a composite knot must either permute cyclically its prime summands, or act 
as a symmetry of one prime summand while permuting the remaining ones. In 
particular the quotient knot cannot be trivial.
\qed

\bigskip

The following is a key lemma for the proofs of Theorems \ref{thm:twins} and 
\ref{thm:rotations}. 

\smallskip

\begin{Lemma}\label{lem:companion} 
Let $K$ be a knot admitting a rotational symmetry $\psi$ of order $p>2$ and 
consider the $JSJ$-decomposition of its exterior $E(K) = \S^3 \setminus\UU(K)$.

\item{(i)} $T$ is a torus of the decomposition which does not separate
$\partial E(K)$ from $Fix(\psi)$ if and only if the orbit $\psi T$ has $p$ 
elements.

\item{(ii)} Under the assumption that $\psi$ has trivial quotient, each torus 
which separates $\partial E(K)$ from $Fix(\psi)$ corresponds to a prime 
companion of $K$ on which $\psi$ acts with trivial quotient.
\end{Lemma}

\demo 
Let $T$ be a torus of the $JSJ$-decomposition of $E(K)$ considered as a torus 
inside $S^3$: $T$ separates the $3$-sphere into a solid torus containing $K$ 
and the exterior of a non trivial knot $K_T$ which is a companion of $K$. Note 
that, since the order of the symmetry $\psi$ is $>2$, its axis cannot meet $T$. 
Assume that the axis $Fix(\psi)$ of the symmetry is contained in the solid 
torus. 

If the orbit of $T$ under $\psi$ does not contain $p$ elements, then a 
non-trivial power of $\psi$ leaves $T$ invariant, and thus it also leaves the 
solid torus and the knot exterior invariant. The restriction of this power of 
$\psi$ to the solid torus acts as a rotation of order $m >1$ around its core 
and leaves invariant each meridian. This non-trivial power of $\psi$ would then 
be a rotational symmetry about the non trivial knot $K_T$ which is absurd 
because of the proof of the Smith's conjecture (see \cite{MB}). 

For the reverse implication, it suffices to observe that the geometric pieces 
of the decomposition containing $\partial E(K)$ and $Fix(\psi)$ must be 
invariant by $\psi$, and so must be the unique geodesic segment joining the 
corresponding vertices in the tree dual to the decomposition.

For the second part of the Lemma, note that $K_T/\psi$ is a companion of 
$K/\psi$, which is trivial by hypothesis. In particular $K_T/\psi$ is also 
trivial and thus, by Lemma \ref{lem:prime}, must be prime.
\qed

\bigskip

The following lemma gives a weaker version of assertion (ii) of Theorem 
\ref{thm:rotations} under a commutativity hypothesis:

\smallskip

\begin{Lemma}\label{lem:two rotations}
Let $K$ be a prime knot admitting two commuting rotational symmetries $\psi$ 
and $\varphi$ of orders $p,q>2$. Then:

\item{(i)} The fixed-point sets of $\psi$ and $\varphi$ are contained in the
same geometric component of the $JSJ$-decomposition for $E(K)$;

\item{(ii)} If $\psi$ has trivial quotient and $p \not = q$, the fixed-point 
sets of $\psi$ and $\varphi$ sit in the component which contains 
$\partial E(K)$.
\end{Lemma}

\demo

\noindent{\bf Part (i)}  Let $v_{\psi}$ (respectively $v_{\varphi}$) the
vertex of the graph $\Gamma_K$ dual to the $JSJ$-decomposition of $E(K)$ 
corresponding to the geometric component containing $Fix(\psi)$ (respectively 
$Fix(\varphi)$). Since the two rotational symmetries commute, $\psi$ 
(respectively $\varphi$) must leave $Fix(\varphi)$ (respectively $Fix(\psi)$) 
invariant, and so the geodesic segment of $\Gamma_K$ joining $v_{\psi}$ to 
$v_{\varphi}$ must be fixed by the induced actions of $\psi$ and $\varphi$ on 
$\Gamma_K$. If this segment contains an edge $e$, the corresponding $JSJ$-torus
$T$ in $E(K)$ cannot separate both $Fix(\varphi)$ and $Fix(\psi)$ from 
$\partial E(K)$. This would contradict part (i) of Lemma \ref{lem:companion}.

\medskip

\noindent{\bf Part (ii)} Let $M$ be the $p$-fold cyclic branched cover of $K$ 
and let $h$ be the associated covering transformation. According to Remark 
\ref{rem:lift} the lift $\tilde \psi$ of $\psi$ to $M$ is the deck 
transformation of a cyclic cover of ${\S}^3$ branched along a knot $K'$. Note 
that both $\tilde \psi$ and $\tilde \varphi$ (the lift of $\varphi$ to $M$) 
commute on $M$ with the covering transformation $h$. In particular 
$\tilde \varphi$ and $h$ induce commuting rotational symmetries of $K'$ with 
order $q$ and $p$ respectively. According to part $(i)$, $Fix(\varphi)$ and 
$Fix(h)$ belong to the same piece of the $JSJ$-decomposition of $M$. Since 
$Fix(h)$ maps to $K$ and $p \not = q$, $Fix(\varphi)$ sits in the $JSJ$-piece 
of $E(K)$ which contains $\partial E(K)$ and the conclusion follows since 
$Fix(\psi)$ belongs to the same $JSJ$-piece as $Fix(\varphi)$.
\qed

\bigskip

\begin{Lemma}\label{lem:torus} 
Let $K$ be a knot admitting a rotational symmetry $\psi$ with trivial quotient 
and of order $p>2$. Let $M$ be the $p$-fold cyclic branched cover of $K$ and 
denote by $\pi:M \longrightarrow(\S^3,K)$ the associated branched cover. Let 
$T$ be a torus in the $JSJ$-collection of tori of $E(K)$.

\item{(i)} The torus $T$ is left invariant by $\psi$ if and only if 
$\pi^{-1}(T)$ is connected.

\item{(ii)} If $\pi^{-1}(T)$ is connected, then the companion $K_T$ of $K$ 
corresponding to $T$ is prime and the winding number of $T$ with respect to $K$ 
is prime with $p$, so in particular it is not zero.

\item{(iii)} The torus $T$ is not left invariant by $\psi$ if and only if 
$\pi^{-1}(T)$ has $p$ components.
\end{Lemma}

\demo

\noindent{\bf Part (i)}. According to Remark \ref{rem:lift}, the $p$-fold 
cyclic branched cover $M$ of $K$ admits two commuting diffeomorphisms of order 
$p$, $h$ and $h'=\tilde\psi $, such that: $(M,Fix(h))/<h> \cong (\S^3,K)$ on 
which $h'$ induces the $p$-rotational symmetry $\psi$ with trivial quotient, 
and $(M,Fix(h'))/<h'> \cong (\S^3,K')$ on which $h$ induces a $p$-rotational 
symmetry $\psi'$ with trivial quotient. The preimage $\pi^{-1}(T) = \tilde T$ 
is connected if and only if it corresponds to a torus $\tilde T$ of the 
$JSJ$-decomposition of $M$ which is left invariant by $h$. If, by 
contradiction, $\psi$ does not leave $T$ invariant, then the $h'$-orbit of 
$\tilde T$ consists of $m>1$ elements. Cutting $M$ along these $m$ separating 
tori, one gets $m+1$ connected components. 

\smallskip

\begin{Claim}\label{claim:component} 
Both $Fix(h)$ and $Fix(h')$ must be contained in the same connected component.
\end{Claim}

\demo
The diffeomorphism $h'$ cyclically permutes the $m$ connected components which 
do not contain $Fix(h')$. Since $h$ and $h'$ commute, $h$ leaves invariant each 
of these $m$ components and it acts in the same way on each of them (that is, 
the restrictions of $h$ to each component are conjugate). Since the set 
$Fix(h)$ is connected, the claim follows.
\qed

\bigskip

The $m$ components permuted by $h'$ project to a connected submanifold of the 
exterior $E(K')$ of the knot $K'$ with connected boundary the image $T'$ of $\tilde T$. This submanifold  is invariant by the action of $\psi'$ 
but does not contain $Fix(\psi')$. This contradicts Lemma 
\ref{lem:companion}(i). To conclude the proof of Lemma \ref{lem:torus} (i), it 
suffices to observe that $h$ and $h'$ play symmetric roles.

\medskip

\noindent{\bf Part (ii)} The first part of assertion (ii) is a straightforward 
consequence of assertion (i) and of Lemma \ref{lem:companion}. The second part 
follows from the fact that for $\pi^{-1}(T)$ to be connected, the winding 
number of $T$ and $p$ must be coprime.

\medskip

\noindent{\bf Part (iii)} is a  consequence of the proof of part (i) of Lemma 
\ref{lem:companion} and of the fact that $h$ and $h'$ play symmetric roles.
\qed

\bigskip

{\bf Proof of Theorem \ref{thm:rotations}.} The proof is achieved in three 
steps.

\medskip

\noindent {\bf Step 1.}  Theorem \ref{thm:rotations} is true under the 
assumption that the rotational symmetries commute pairwise.

In this case, assertion (i) is the statement of Proposition \ref{prop:commuting 
rotations}. Assertion (ii) follows from Lemma \ref{lem:two rotations}. 

\bigskip

\noindent{\bf Step 2.} Theorem \ref{thm:rotations} is true under the assumption 
that every companion of $K$ is prime (i.e. $K$ is \emph{totally prime}) and has 
non vanishing winding number (i.e. $K$ is \emph{pedigreed}).

Assume that we are in the hypotheses of Theorem \ref{thm:rotations}. Then Lemma 
\ref{lem:prime} assures that $K$ is a prime knot. If $K$ is also totally prime  
and pedigreed then M. Sakuma \cite[Thm4 and Lemma 2.3]{Sa2} proved that, up to 
conjugacy, the rotational symmetries belong either to a finite cyclic subgroup 
or to an $S^1$-action in $Diff^{+,+}(S^3,K)$. Thus after conjugacy, step 1 
applies. For part (ii) note that the distances of the fixed point set of the 
symmetries to the vertex containing $\partial E(K)$ in the $JSJ$-graph 
$\Gamma_ K$ do not change by conjugacy.

\bigskip

\noindent {\bf Step 3.} Reduction of the proof to step 2.

If $K$ is not totally prime or pedigreed, then it is non-trivial. We shall
construct a non trivial, totally prime and pedigreed knot verifying the 
hypothesis of Theorem \ref{thm:rotations}. Assertion (i) then follows by 
contradiction. For Assertion (ii) we need to verify that the construction does 
not change the distance of the pieces containing the axes of rotations to the 
root containing $\partial E(K)$. Roughly speaking we consider the $JSJ$-tori 
closest to $\partial E(K)$ and corresponding either to non-prime or to winding 
number zero companions. Then we cut $E(K)$ along these tori and keep the 
component $W$ containing $\partial E(K)$ and suitably Dehn-fill $W$ along these 
tori to get the exterior of a non-trivial knot $\hat K$ in $\S^3$, which 
verifies Sakuma's property. 

More precisely, let $\Gamma_K$ be the tree dual to the $JSJ$-decomposition of 
$E(K)$ and let $\Gamma_0$ be its maximal (connected) subtree with the following 
properties: 
\begin{itemize}

\item $\Gamma_0$ contains the vertex $v_\partial$ corresponding to the 
geometric piece whose boundary contains $\partial E(K)$. Note that the
geometric piece of the decomposition corresponding to $v_\partial$ cannot be a
composing space for $K$ is prime;

\item no vertex of $\Gamma_0$ corresponds to a composing space (i.e. a space 
homeomorphic to a product $S^1 \times B$ where $B$ is an $n$-punctured disc 
with $n \geq 2$);

\item no edge of $\Gamma_0$ corresponds to a torus whose meridian has linking 
number $0$ with $K$.

\end{itemize}
Denote by $X(\Gamma_0)$ the submanifold of $E(K)$ corresponding to $\Gamma_0$.

The following claim describes certain properties of $X(\Gamma_0)$ with respect 
to a rotational symmetry $\psi$ of $(\S^3,K)$.

\begin{Claim}\label{claim:sym}  
Let $\psi$ be a rotational symmetry of $(\S^3,K)$ with order $p >2$ and trivial
quotient. Then:

\item{(i)} The fixed-point set of $\psi$ is contained in $X(\Gamma_0)$.

\item{(ii)} The tree $\Gamma_0$ is invariant by the automorphism of $\Gamma_K$ 
induced by $\psi$ and the submanifold $X(\Gamma_0)$ is invariant by $\psi$.
\end{Claim}

\demo

\noindent{\bf Assertion {(i)}.} Let $\gamma$ be the unique geodesic segment in 
$\Gamma_K$ which joins the vertex $v_\partial$ to the vertex corresponding to 
the geometric piece containing $Fix(\psi)$ (see Lemma \ref{lem:JSJ}; note that 
here we use $p>2$). According to assertion (ii) of Lemma \ref{lem:companion}, 
no vertex along $\gamma_i$ can be a composing space. Since the linking number 
of $K$ and $Fix(\psi)$ must be coprime with $p$, no torus corresponding to an 
edge of $\gamma$ can have winding number $0$ (see Lemma \ref{lem:torus}). 

\medskip

\noindent{\bf Assertion {(ii)}.} This is just a consequence of the maximality 
of $\Gamma_0$ and the fact that elements of the group $\langle\psi \rangle$ 
generated by $\psi$ must preserve the $JSJ$-decomposition of $E(K)$ and the 
winding numbers of the $JSJ$-tori, as well as send composing spaces to 
composing spaces.
\qed

\bigskip

Let $\pi :M_{p}(K) \longrightarrow(\S^3,K)$ be the $p$-fold cyclic 
branched cover. Let $T$ be a torus of the $JSJ$-collection of tori for $E(K)$. 
Denote by $E_T$ the manifold obtained as follows: cut $E(K)$ along $T$ and 
choose the connected component whose boundary consists only of $T$. Note that 
$E_T$ is the exterior of the companion $K_T$ of $K$ corresponding to $T$.

\smallskip

\begin{Claim}\label{claim:meridian-longitude}
Let $T$ be a torus of $\partial X(\Gamma_0)\setminus\partial E(K)$. The 
preimage $\pi^{-1}(T)$ consists of $p$ components, each bounding a copy of 
$E_T$ in $M_{p}(K)$. In particular, there is a well-defined
meridian-longitude system $(\mu_T,\lambda_T)$ on each boundary component of
$X(\Gamma_0)$, different from $\partial E(K)$, which is preserved by taking the
$p$-fold cyclic branched covers.
\end{Claim}
\demo
According to Lemma \ref{lem:torus}, the preimage of $T$ is either connected or 
consists of $p$ components. If the preimage of $T$ were connected, the tree 
$\Gamma_0$ would not be maximal according to Lemma \ref{lem:torus}(ii). The 
remaining part of the Claim is then easy.
\qed

\bigskip

We wish now to perform Dehn fillings on the boundary of $X(\Gamma_0)$ in order 
to obtain a totally prime and pedigreed knot admitting pairwise distinct 
rotational symmetries with trivial quotients. On each component $T$ of
$\partial X(\Gamma_0)\setminus\partial E(K)$ we fix the curve $\alpha_n = 
\lambda_T+n\mu_T$.

\smallskip

\begin{Claim}\label{claim:surgery} 
For all but finitely many $n \in \Z$ the Dehn filling of each component $T$ of 
$\partial X(\Gamma_0)\setminus\partial E(K)$ along the curve $\alpha_n$
produces the exterior of a non-trivial, prime and pedigreed knot $\hat K$ in 
$\S^3$.
\end{Claim}

\demo
Note that by the choice of surgery curves the resulting manifold 
$\hat X(\Gamma_0)$ is the exterior of a knot $\hat K$ in the $3$-sphere, i.e. 
$\hat X(\Gamma_0)\subset \S^3$, and thus is irreducible. We distinguish two 
cases:

\medskip

\noindent{\bf {(1)}} The $JSJ$-component $X_T$ of $X(\Gamma_0)$ adjacent to $T$ 
is Seifert fibred. Then, by the choice of $\Gamma_0$, $X_T$ is a cable space 
(i.e. the exterior of a $(a,b)$-torus knot in the solid torus bounded by $T$ in
$\S^3$). Moreover the fiber $f$ of the Seifert fibration of $X(\Gamma_0)$ is 
homologous to $a\mu_{T} + b\lambda_{T}$ on $T$ and the intersection number
$\vert \Delta(f,\mu_T) \vert = b > 1$. The intersection number of the filling 
curve $\alpha_n$ with the fiber $f$ is then $\vert \Delta(f,\alpha_n) \vert =
\vert na -b \vert$ and is $> 1$  for all but finitely many $n \in \Z$. In this 
case the resulting manifold $X_T(\alpha_n)$ is the exterior of a non trivial 
torus knot which is prime and pedigreed \cite{CGLS}.

\medskip

\noindent {\bf{( 2)}} The $JSJ$-component $X_T$ of $X(\Gamma_0)$ adjacent to 
$T$ is hyperbolic. By Thurston's hyperbolic Dehn filling theorem 
\cite[Chap. 5]{Th1} (see also \cite[Appendix B]{BoP}) for all but finitely many 
$n \in \Z$ the Dehn filling of each component $T \subset \partial X_T \cap 
(\partial X(\Gamma_0)\setminus\partial E(K))$ along the curve $\alpha_n$ 
produces a hyperbolic manifold $X_T(\alpha_n)$ with finite volume.

\medskip

Therefore for all but finitely many $n$'s $\in \Z$ the Dehn filling of each 
component $T \subset \partial X(\Gamma_0)\setminus\partial E(K)$ along the 
curve $\alpha_n$  produces a $\partial$-irreducible $3$-manifold $\hat 
X(\Gamma_0) \subset \S^3$ such that each Seifert piece of its 
$JSJ$-decomposition is either a Seifert piece of $X(\Gamma_0)$ or a non-trivial 
torus knot exterior. Hence it corresponds to the exterior of a non-trivial knot 
$\hat{K} \subset \S^3$ which is totally prime. It is also pedigreed by the 
choice of $\Gamma_0$.
\qed

\bigskip

Let $\psi$ a rotational symmetry of $(\S^3,K)$ with order $p >2$. Then the 
restriction ${\psi}_{\vert_{X(\Gamma_0)}}$, given by Claim \ref{claim:sym} 
extends to $\hat X(\Gamma_0)$, giving a $p$-rotational symmetry $\hat\psi$ of 
the non-trivial, totally prime and pedigreed knot $(\S^3,\hat K)$. In order to 
be able to apply  step 2 to the knot $\hat K$ and the induced rotational 
symmetries, we still need to check that the rotational symmetry $\hat \psi$ has 
trivial quotient when $\psi$ has trivial quotient. This is the aim of the 
following:

\smallskip

\begin{Claim}\label{claim:quotient} 
If the knot $K/\psi$ is trivial, then the knot $\hat K/\hat\psi$ is trivial.
\end{Claim}

\demo 
Let $\pi:M_{p}(K) \longrightarrow(\S^3,K)$ be the $p$-fold cyclic 
branched cover. Let $h$ be the deck transformation of this cover and $h'$ 
the lift of $\psi$. According to Remark \ref{rem:lift}, $h'$ is the deck 
transformation for the $p$-fold cyclic cover of the $3$-sphere branched along 
a knot $K'$. Note that, by Claim \ref{claim:meridian-longitude},
$M_{p}(K)\setminus\pi^{-1}(X(\Gamma_0) \cup \UU(K))$ is a disjoint union of $p$ 
copies of $E(K) \setminus X(\Gamma_0)$. It follows that the $p$-fold cyclic 
branched cover $M_{p}(\hat K)$ of $\hat K$ is the manifold obtained by a 
$(\lambda_T+n\mu_T)$-Dehn filling on all the boundary components of 
$\pi^{-1}(X(\Gamma_0) \cup \UU(K))$. The choice of the surgery shows that both 
$h$ and $h'$ extend to diffeomorphisms $\hat h$ and $\hat h'$ of order 
$p$ of $M_{p}(\hat K)$. By construction we have that $M_{p}(\hat K)/<\hat
h> \cong \S^3$. In the same way $M_{p}(\hat K)/<\hat h'>$ is obtained
from $M_{p}(K)/<h'> \cong\S^3$ by cutting off a copy of $E(K) \setminus
X(\Gamma_0)$ and Dehn filling along $\partial X(\Gamma_0)$. The choice of
the surgery curve assures that the resulting manifold is again $\S^3$ and the 
conclusion follows from Remark \ref{rem:lift}.
\qed

\bigskip

From the non-trivial prime knot $K$, we have thus constructed a non-trivial, 
totally prime and pedigreed knot $\hat K$ which has the property that every 
rotational symmetry $\psi$ of $K$ with trivial quotient and order $>2$  
induces a rotational symmetry $\hat \psi$ of $\hat K$ with trivial quotient and 
the same order. Moreover by the choice of the Dehn filling curve in the 
construction of $\hat K$, the vertex containing $Fix(\hat \psi)$ remains at the 
same distance from the vertex containing $\partial E(\hat K)$ in the 
$JSJ$-tree $\Gamma_ {\hat K}$ as the vertex containing $Fix(\psi)$ from the 
vertex containing $\partial E(K)$ in the $JSJ$-tree $\Gamma_ K$. Then the 
conclusion is a consequence of step 2.\qed

\section{Twins of a prime knot}

In this section we prove Theorem \ref{thm:twins}. If $K$ is trivial, the
theorem is a consequence of the proof of Smith's conjecture (see 
\cite{MB}). We shall thus assume in the remaining of this section that $K$ is 
non trivial and $p$ is an odd prime number.

Let $M$ be the common $p$-fold cyclic branched cover of two prime knots $K$ and 
$K'$ in $\S^3$. Let $h$ and $h'$ be the deck transformations for the coverings 
of $K$ and $K'$ respectively. By the orbifold theorem \cite{BoP}, see also 
\cite{CHK} one can assume that $h$ and $h'$ act \emph{geometrically} on the 
geometric pieces of the $JSJ$-decomposition of $M$, i.e. by isometries on the 
hyperbolic pieces and respecting the fibration on the Seifert fibred ones.

The following lemma describes the Seifert fibred pieces of the 
$JSJ$-decomposition of the $p$-fold branched cyclic cover $M$ (see also
\cite{Ja} and \cite[Lemma 2]{Ko}).

\smallskip

\begin{Lemma}\label{lem:seifert}
Let $p$ be an odd prime integer and let $M$ be the $p$-fold cyclic branched 
cover of $\S^3$ branched along a prime, satellite knot $K$. If $V$ is a Seifert 
piece in the $JSJ$-decomposition for $M$. Then the base $B$ of $V$ can be:

\begin{enumerate}

\item A disc with $2$, $p$ or $p+1$ singular fibres;

\item A disc with $1$ hole, i.e. an annulus, with $1$ or $p$ singular fibres;

\item A disc with $p-1$ holes and $1$ singular fibre;

\item A disc with $p$ holes and $1$ singular fibre;

\item A disc with $n$ holes, $n\ge2$.

\end{enumerate}
\end{Lemma}

\demo
It suffices to observe that $V$ projects to a Seifert fibred piece $V'$ of the 
Bonahon-Siebenmann decomposition for the orbifold $\OO_p(K)$. There are four 
possible cases:

\noindent{\bf (a)} $V'$ contains $K$: $V'$ is topologically a non trivially 
fibred solid torus and $K$ is a regular fibre of the fibration, i.e. a torus 
knot $K(a,b)$, since it cannot be the core of the fibred solid torus. The knot
$K$ lifts to a singular fibre of order $p$ if $p$ does not divide $ab$ and to 
a regular fibre otherwise. The core of the solid torus is a singular fibre of 
order -say- $a$. It lifts to a regular fibre if $a=p$, a singular fibre of 
order $a/(a,p)$ if $p$ does not divide $b$, or to $p$ singular fibres of order 
$a$ if $p$ divides $b$. Thus $V$ has $p$ boundary components if $p$ divides $a$ 
and $1$ otherwise. An Euler characteristic calculation shows that $B$ is either 
a disc with $2$ or $p$ singular fibres, or a disc with $p-1$ holes and with at 
most $1$ singular fibre.

\noindent{\bf (b)} $V'$ is the complement of a torus knot $K(a,b)$ in $\S^3$. 
In this case, $V$ is either a copy of $V'$, and $B$ is a disc with $2$ singular
fibres or $V$ is a true $p$-fold cover of $V'$. In this case $V$ has exactly
one boundary component. Reasoning as in case (a), we see that the two singular
fibres of $V'$ must lift to either $2$ singular fibres, or $1$ regular fibre
and $p$ singular fibres or $1$ singular fibre and $p$ singular fibres. In
particular $B$ is a disc with $2$, $p$ or $p+1$ singular fibres.

\noindent{\bf (c)}  $V'$ is the complement of a torus knot $K(a,b)$ in a solid 
torus, i.e. a cable space, and its base is an annulus with $1$ singular fibre. 
Reasoning as in (b) we find that $B$ can be a disc with $1$ hole and $1$ or $p$ 
singular fibres or a disc with $p$ holes and at most $1$ singular fibre.

\noindent{\bf (d)}  $V'$ is a composing space with at least $3$ boundary
components and thus so is $V$. More precisely, note that either $V'$ lifts to
$p$ disjoint copies of itself, or $V$ and $V'$ are homeomorphic and $V'$ is 
obtained by quotienting $V$ via the $p$-translation along the $\S^1$ fibre. In 
this case $B$ is a disc with at least $2$ holes.

This analysis ends the proof of Lemma \ref{lem:seifert}.
\qed

\bigskip

\begin{Proposition}\label{prop:subtree}
Let $M$ be the common $p$-fold cyclic branched cover of two prime knots $K$ and 
$K'$ in $\S^3$, $p$ an odd prime number, and let $h$ be the deck transformation 
for the covering of $K$. Let $\Gamma$ be the tree dual to the
$JSJ$-decomposition of $M$. The deck transformation $h'$ for the covering of
$K'$ can be chosen (up to conjugacy) in such a way that:

\item{(i)} There exists a subtree $\Gamma_f$ of $\Gamma$ on which the actions
induced by $h$ and $h'$ are trivial;

\item{(ii)} The vertices of $\Gamma$ corresponding to the geometric pieces of
the decomposition which contain $Fix(h)$ and $Fix(h')$ belong to $\Gamma_f$;

\item{(iii)} Let $M_f$ the submanifold of $M$ corresponding to $\Gamma_f$. The
restrictions of $h$ and $h'$ to $M_f$ commute.
\end{Proposition}

\demo 
The proof relies on the study of the actions of the two covering 
transformations $h$ and $h'$ on the $JSJ$-decomposition of the common $p$-fold  
cyclic branched covering $M$. Since $\Gamma$ is finite, the group generated by 
the tree automorphisms induced by $h$ and $h'$ is finite as well. Standard 
theory of group actions on trees assures that a finite group acting on a tree 
without inversion must have a global fixed point and that its fixed-point set 
is connected. Thus part (i) of the proposition follows, using the fact that $h$ 
and $h'$ have odd orders.

\medskip

Choose now $h'$, up to conjugacy in $Diff^+(M)$, in such a way that
$\Gamma_f$ is maximal. We want to show that, in this case, $M_f$ contains 
$Fix(h)$ and $Fix(h')$. Assume by contradiction that the vertex $v_h$ of 
$\Gamma$ corresponding to the geometric piece containing $Fix(h)$, whose
existence is ensured by Lemma \ref{lem:JSJ}, does not belong to $\Gamma_f$. Let 
$\gamma_h$ the unique geodesic path in $\Gamma$ connecting $v_h$ to $\Gamma_f$. 
Let $e_h$ the edge in $\gamma_h$ adjacent to $\Gamma_f$ and denote by $T$ the 
corresponding torus of the $JSJ$-collection of tori for $M$. Let $U$ be the 
connected component of $M\setminus T$ which contains $Fix(h)$. Consider the 
$\langle h,h'\rangle$-orbit of $U$. This orbit is the disjoint union of $h$ 
(and $h'$) orbits of $U$. Remark that the $h$-orbit of $U$ is $\{U\}$.

\smallskip

\begin{Claim}\label{claim:orbit}
The orbit $\langle h,h'\rangle U$ must contain an $h$-orbit, different from 
$\{U\}$ and containing a unique element.
\end{Claim}

\demo 
Otherwise all the $h$-orbits in $\langle h,h'\rangle U$ different from $\{U\}$ 
would have $p$ elements, since $p$ is prime. In particular, the cardinality of 
$\langle h,h'\rangle U$ would be of the form $kp+1$. This implies that at least 
one of the $h'$-orbits in $\langle h,h'\rangle U$ must contain one single 
element $U'$. Up to conjugacy with an element of $\langle h,h'\rangle$ (whose 
induced action on $\Gamma_f$ is trivial), we can assume that $U=U'$, 
contradicting the hypothesis that $h'$ was chosen up to conjugacy in such a 
way that $\Gamma_f$ is maximal.
\qed

\bigskip

Let $U'\neq U$ the element of $\langle h,h'\rangle U$ such that $h(U')=U'$. 
Note that $U$ and $U'$ are homeomorphic since they belong to the same 
$\langle h,h'\rangle$-orbit.

\smallskip

\begin{Claim}\label{claim:knot complement}
$U$ is homeomorphic to the exterior $E({\KK})$ of a knot $\KK \subset \S^3$
admitting a free symmetry of order $p$.
\end{Claim}

\demo 
The first part of the Claim follows from the fact that, by maximality of 
$\Gamma_f$, $h'$ cannot leave $U$ invariant, so must freely permute $p$ copies
of $U$ belonging to $\langle h,h'\rangle U$. Thus $U$ must appear as a union of
geometric pieces of the $JSJ$-splitting of $E(K')$. The second part follows 
from the fact that $h$ must act freely on $U'$ which is homeomorphic to $U$.
\qed

\bigskip

\begin{Remark}\label{rem:free quotient} 
Note that the quotient of $U$ by the action of its free symmetry of order $p$ 
is also a knot exterior because $h$ acts freely on $U'$ and $U'$ must project 
to a union of geometric pieces of the $JSJ$-splitting of $E(K)$.
\end{Remark}

\medskip

\begin{Claim}\label{claim:non-free quotient}
$U$ admits a rotational symmetry of order $p$ whose quotient $U/\langle 
h\rangle$ is topologically a solid torus.
\end{Claim}

\demo
The quotient $U/\langle h\rangle$ is obtained by cutting $\S^3$ along an 
essential torus in $E(K)$. Since $K \subset U/\langle h\rangle$, it must be a 
solid torus. 
\qed

\bigskip

It follows from Claim \ref{claim:non-free quotient} and Lemma \ref{lem:prime} 
that the knot $\KK$ is prime. Moreover, according to Claims \ref{claim:knot 
complement} and \ref{claim:non-free quotient}, $\KK$ admits a rotational 
symmetry and a free symmetry, both of order $p$. This is however impossible  
because M. Sakuma \cite[Thm. 3]{Sa2} showed that a prime knot can only have one 
symmetry of odd order up to conjugacy. This contradiction proves part (ii) of 
Proposition \ref{prop:subtree}.

\medskip

To prove part (iii) we shall consider two cases, according to the structure
of $\Gamma_f$.

\medskip

\noindent {\bf Case (a)}: {\it $\Gamma_f$ contains an edge.}
Choose an edge in $\Gamma_f$ and let $T$ be the corresponding torus in the
$JSJ$-collection of tori for $M$. Let $V$ be a geometric piece of the 
$JSJ$-decomposition of $M$ adjacent to $T$. Then Lemma \ref{lem:commutation} 
below together with a simple induction argument show that $h'$ can be chosen 
(up to conjugacy) in such a way that its restriction to $M_f$ commutes with the 
restriction of $h$.

\smallskip

\begin{Lemma}\label{lem:commutation}
If the covering transformations $h$ and $h'$ preserve a $JSJ$-torus $T$ of $M$ 
then, up to conjugacy in $Diff^{+}(M)$, $h$ and $h'$ commute on the union of 
the geometric components of the $JSJ$-decomposition adjacent to $T$.
\end{Lemma}

\demo
First we show that $h$ and $h'$ commute on each geometric component adjacent to 
$T$. Since $h$ and $h'$ preserve the orientation of $M$, we deduce that 
$h(V)=V$ and $h'(V)=V$, and that $h$ and $h'$ act geometrically on the 
geometric piece $V$. A product structure on $T$ can always be induced by the 
geometric structure on $V$: either by considering the induced Seifert fibration 
on $T$ if $V$ is Seifert fibred, or by identifying $T$ with a section of a cusp 
in the complete hyperbolic manifold $V$. Since $h$ and $h'$ are isometries of 
order $p$, for such a product structure on $T$ they act as (rational) 
translations, i.e. their action on $T={\S}^1\times{\S}^1$ is of the form 
$(\zeta_1,\zeta_2) \mapsto (e^{2i\pi r_1/p}\zeta_1,e^{2i\pi r_2/p}\zeta_2)$, 
where $p$ and at least one between $r_1$ and $r_2$ are coprime. Thus $h$ and 
$h'$ commute on $T$.

\medskip

If $V$ is hyperbolic, we have just seen that $h$ and $h'$ are two isometries of 
$V$ which commute on the cusp corresponding to $T$. Thus they must commute on 
$V$.

\medskip

If $V$ is Seifert fibred, then the Seifert fibration is unique up to isotopy, 
and $h$ and $h'$ preserve this fibration. 

\smallskip

\begin{Remark}
Note that the quotient of $V$ by a fiber-preserving diffeomorphism of
finite order $h$ only depends on the combinatorial behaviour of $h$, i.e. its
translation action along the fibre and the induced permutation on cone points 
and boundary components of the base. In particular, the conjugacy class of $h$ 
only depends on these combinatorial data. Note moreover that two geometric
symmetries having the same combinatorial data are conjugate via a
diffeomorphism isotopic to the identity.
\end{Remark}

\medskip

Since the translation along the fibres commutes with every fiber-preserving 
diffeomorphism of $V$, it suffices to see whether $h$ and $h'$ commute, up to a 
conjugation of $h'$, on the base $B$ of $V$. It is enough then to consider the 
possible actions of order $p$ on the possible bases. According to Lemma 
\ref{lem:seifert} the possible actions of $h$ and $h'$ are described below:

\begin{enumerate}

\item If $B$ is a disc with $2$ singular fibres, or an annulus with $1$ 
singular fibre, or a disc with $n$ holes, $n\neq p$, or a disc with $p-1$ holes
and $1$ singular fibre, then the action on $B$ is necessarily trivial and there 
is nothing to prove. Note that, according to the proof of Lemma 
\ref{lem:seifert}, if $B$ is a disc with $p-1$ holes with one singular fibre, 
no boundary torus is left invariant, so this possibility in fact does not 
occur.

\item If $B$ is a disc with $p$ holes and $1$ singular fibre or a disc with 
$p+1$ singular fibres, then the only possible action is a rotation about a 
singular fibre cyclically permuting the holes or the remaining singular fibres. 

\item If $B$ is a disc with $p$ singular fibres then the action must be a
rotation about a regular fibre which cyclically exchanges the singular fibres.

\item If $B$ is an annulus with $p$ singular fibres the action must be a free
rotation cyclically exchanging the singular fibres. Note that in the three
latter cases the action can never be trivial on the base.

\item If $B$ is a disc with $n$ holes then two situations can arise: either the
action is trivial on the base (case (d) in the proof of Lemma 
\ref{lem:seifert}; note that in case (a), when $n=p-1$, all boundary components 
must be cyclically permuted), or $n=p$ and the action is a rotation about a 
regular fibre which cyclically permutes the $p$ holes (see part (c) of Lemma 
\ref{lem:seifert}).

\end{enumerate}

We shall now show that, if both $h$ and $h'$ induce non trivial actions on the
base of $V$, then, up to conjugacy, $h$ and $h'$ can be chosen so that their 
actions on $B$ coincide. Note that for $h$ and $h'$ to commute it suffices that 
the action of $h'$ on $B$ coincides with the action of some power of $h$, 
however this stronger version will be needed in the proof of Corollary
\ref{cor:extension}.

First of all remark that, if $B$ is a disc with $p+1$ singular fibres (case 2) 
and $h$ and $h'$ leave invariant distinct singular fibres, then all the
singular fibres must have the same order (in fact, must have the same
invariants). This means that, after conjugating $h'$ by a homeomorphism of $V$
which is either an isotopy exchanging two regular fibres or a Dehn twist along
an incompressible torus exchanging two singular fibres, one can assume that, in
cases 2 and 3, $h$ and $h'$ leave set-wise invariant the same fibre. Note that
this homeomorphism is isotopic to the identity on $\partial V$ and thus extends
to $M$. In fact, using Lemma \ref{lem:seifert} one can show that the fibres
cannot all have the same order.

Since the actions of $h$ and $h'$ consist in permuting exactly $p$ holes or 
singular fibres, it suffices to conjugate $h'$ via a homeomorphism of $V$ 
(which is a composition of Dehn twists along incompressible tori) in such a way 
as to exchange the order of the holes or singular fibres so that $h'$ and $h$
cyclically permute them in the same order. Note that in the case of singular 
fibres this product of Dehn twists is isotopic to the identity on $\partial V$ 
and thus extends to $M$. In the case of holes, the product of Dehn twists 
extends to $M$ since it induces the identity on the fundamental groups of the 
tori of $\partial V$ and the connected components of $M\setminus V$ adjacent to 
boundary tori different from $T$ are necessarily homeomorphic. 

Once the two diffeomorphisms $h$ and $h'$ commute on the two geometric pieces 
adjacent to $T$, the commutation can be extended on a product neighborhood of 
$T$, since the two finite abelian groups generated by the restrictions of $h$ 
and $h'$ on each side of $T$ have the same action on $T$. Indeed, the slope of 
the translation induced by $h'$ on $T$ has been left unchanged by the 
conjugation.
\qed

\bigskip

\begin{Remark}\label{rem:same action}
Note that in case 1 of the proof of the above Lemma, the actions of $h$ and
$h'$ must coincide after taking a power, i.e. $h$ and $h'$ generate the same
cyclic group. This is not necessarily true in the remaining cases, even if $h$ 
and $h'$ induce the same action on $B$. Indeed, they can induce different
translations along the fibres. Nevertheless, in both cases, to assure that the
actions of $h$ and $h'$ coincide on $V$, it suffices to check that they 
coincide on $T$.
\end{Remark}

\medskip

\noindent {\bf Case (b)}: {\it $\Gamma_f$ is a single vertex.}
Let $V=M_f$ be the geometric piece corresponding to the unique vertex of 
$\Gamma_f$. If $V=M$, then the result is already known. We can thus assume that
$V\neq M$. According to part (ii) of Proposition \ref{prop:subtree}, we 
can assume that the fixed-point sets of $h$ and $h'$ are contained in $V$.
If $V$ is Seifert fibred then, case (a) of the proof of Lemma \ref{lem:seifert} 
shows that the base $B$ of $V$ is either a disc with $2$ or $p+1$ singular 
fibres, or a disc with $p-1$ holes and with $1$ or $2$ singular fibres. In the 
first case the boundary torus of $V$ is preserved by $h$ and $h'$ and the 
assertion follows from Lemma \ref{lem:commutation}. In the second case the 
action on the base is necessarily a rotation fixing two points (either the 
unique singular fibre and a regular one, or the two singular fibres) and 
cyclically permuting the $p$ boundary components. Then conjugating $h'$ by a 
product of Dehn twists along incompressible tori, which extends to $M$ as in 
the proof of Lemma \ref{lem:commutation}, leads to the desired conclusion.

\medskip

The case where $V$ is hyperbolic is due to B. Zimmermann \cite{Zim1}. We give 
the argument for completeness. Since $V$ is hyperbolic, we consider the group 
$\II_V$ of isometries of $V$ induced by diffeomorphisms of $M$ which leave $V$ 
invariant. Let $\SS$ be the $p$-Sylow subgroup of $\II_V$. Up to conjugacy, we 
can assume that both $h=h_{\vert_V}$ and $h'=h'_{\vert_V}$ belong to $\SS$. If 
the groups $\langle h\rangle$ and $\langle h'\rangle$ generated by $h$ and $h'$ 
are conjugate, we can assume that $h=h'$ and we are done. So we assume that 
$\langle h\rangle$ and $\langle h'\rangle$ are not conjugate. Then it suffices 
to prove that $h'$ normalises $\langle h\rangle$ because each element 
normalising $\langle h\rangle$ must leave invariant $Fix(h)$ and the subgroup 
of $\II_V$ which leaves invariant a simple closed geodesic, like $Fix(h)$, must 
be a finite subgroup of $\Z/2\Z\ltimes(\Q/\Z\oplus\Q/\Z)$. In particular, 
elements of odd order must commute. Assuming that $\langle h\rangle$ and 
$\langle h'\rangle$ are not conjugate, we have that $\langle h\rangle 
\subsetneq \SS$ and, by \cite[Ch 2, 1.5]{Su}, either $\langle h\rangle$ is 
normal in $\SS$ and we have reached the desired conclusion, or there exist an 
element $\hat{h}=ghg^{-1}$, conjugate to $h$ in $\SS$, which normalises 
$\langle h\rangle$ and such that 
$\langle h\rangle \cap\langle{\hat{h}}\rangle=\{1\}$.

We want to show that $h'$ normalises $\langle h\rangle$. Assume, by 
contradiction that $h'$ is not contained in $\langle h,\hat{h}\rangle =
\Z/p\Z\oplus\Z/p\Z$. Then this group is smaller than $\SS$ and again we are 
able to find a new cyclic group $H$ of order $p$ whose intersection with 
$\langle h,\hat{h}\rangle$ is reduced to the identity and which normalises 
$\langle h,\hat{h}\rangle$. Since the order of $H$ is an odd prime number and 
since $\langle h\rangle$ and $\langle\hat{h}\rangle$ are the only subgroups of 
$\langle h,\hat{h}\rangle$ which fix point-wise a geodesic by \cite[Proposition 
4]{MZ}, $H$ would commute with $\langle h,\hat{h}\rangle$ which is a 
contradiction to the structure of a group leaving a geodesic invariant. This 
final contradiction shows that, up to conjugacy, the subgroups $\langle 
h\rangle$ and $\langle h'\rangle$ either commute or coincide on $V$. This 
finishes the proof of Proposition \ref{prop:subtree}.
\qed

\bigskip

The following proposition shows that a prime knot $K$ having a $p$-twin either 
admits a rotational symmetry of order $p$, or a well-specified submanifold 
$E_p(K)$ built up of geometric pieces of the $JSJ$-decomposition of $E(K)$ 
admits a symmetry of order $p$ with non-empty fixed-point set.

\smallskip

\begin{Definition} 
Let $K$ be a prime knot in $\S^3$. For each odd prime number $p$ we define 
$E_p(K)$ to be the connected submanifold of $E(K)$ containing $\partial E(K)$ 
and such that $\partial E_p(K) \setminus \partial E(K)$ is the union of the 
$JSJ$-tori of $E(K)$ with winding number $p$ which are closest to $\partial 
E(K)$.
\end{Definition}

\medskip

\begin{Proposition}\label{prop:orbifold}
Let $K$ be a prime knot and let $p$ be an odd prime number. Then for any 
$p$-twin $K'$, the deck transformation of the branched cover
$M\longrightarrow({\S}^3,K')$ induces on $E_p(K)$ a symmetry of order $p$, with
non-empty fixed-point set and which extends to $\UU(K)$.
\end{Proposition}

\demo
First we show that the deck transformation of the branched cover
$M\longrightarrow({\S}^3,K')$ associated to a $p$-twin of $K$ induces on 
$E_p(K)$ a symmetry of order $p$.
  
Let $K'$ be a $p$-twin of $K$. Let $h$ and $h'$ be the deck transformations on 
$M$ for the $p$-fold cyclic branched covers of $K$ and $K'$. We shall start by 
understanding the behaviour of $h$ and $h'$ on $M$. We have seen in Proposition
\ref{prop:subtree} that $h$ and $h'$ can be chosen to commute on the 
submanifold $M_f$ of $M$ corresponding to the maximal subtree of $\Gamma$ on 
which both $h$ and $h'$ induce a trivial action. Let $\Gamma_c$ the maximal 
$\langle h,h'\rangle$-invariant subtree of $\Gamma$ containing $\Gamma_f$, such 
that, up to conjugacy, $h$ and $h'$ can be chosen to commute on the 
corresponding submanifold $M_c$ of $M$.

If $M_c = M$ then after conjugation $h'$ commutes with $h$ on $M$, but is
distinct from $h$ because the knots $K$ and $K'$ are not equivalent. Hence it 
induces a rotational symmetry of order $p$ of the pair $(S^3,K)$ and we are 
done.

So we consider now the case where $\partial M_c$ is not empty. It is sufficient 
to show that $E_p(K) \subset M_c/<h>$: then the symmetry of order $p$ induced 
by $h'$ on $M_c/<h>$ must preserve $E_p(K)$ since each $JSJ$-torus of $E(K)$ 
can only be mapped to another torus of the family with the same winding number 
and the same distance from $\partial E(K)$. First we show:

\smallskip

\begin{Claim}\label{claim:permutation}
Let $T$ be a connected component of $\partial M_c$. The $h$-orbit of $T$ 
consists of $p$ elements which are permuted in the same way by $h$ and $h'$.
\end{Claim}

\demo 
Let $T$ be a torus in $\partial M_c$ and let $U$ be the connected component of 
$M\setminus M_c$ adjacent to $T$. Because of Lemma \ref{lem:commutation}, $T$ 
cannot be preserved by both $h$ and $h'$ for else $M_c$ would not be maximal.
Without loss of generality, we can assume that either:
\medskip

\noindent{\bf (a)} $h(T) \neq T$ and $h'(T) \neq T$;

\medskip

\noindent or

\medskip

\noindent{\bf (b)} $h(T)=T$ but $h'(T)\neq T$; in this case since $h$ and
$h'$ commute on $M_c$, we have that $h(h{'}^{\alpha}(U)) = h{'}^{\alpha}(U)$. 
Then part (ii) of Proposition \ref{prop:subtree} implies that $h$ acts freely 
on $h{'}^{\alpha}(U)$ for each $\alpha=0,...,p-1$.

\medskip

In case (a), the orbit of $T$ by the action of the group $\langle h,h'\rangle$ 
consists of $p$ or $p^2$ elements which bound on one side $M_c$ and on the 
other side a manifold homeomorphic to $U$. If the orbit consist of $p$ 
elements, since $h$ and $h'$ commute on $M_c$, up to choosing a different
generator in $\langle h'\rangle$ we can assume that $h$ and $h'$ permute the
elements of the orbit in the same way. Indeed, we have 
$h'h(T)=hh'(T)=h(h^{\alpha}(T))=h^{\alpha}(h(T))$.

If the orbit consist of $p^2$ elements, $U$ is a is a knot exterior and there 
is a well-defined longitude-meridian system on each component of the $\langle 
h,h'\rangle$-orbit of $T$. In particular, there is a unique way to glue a copy 
of $U$ along the projection of $T$ in $M_c/\langle h,h'\rangle$. This implies 
that $h$ and $h'$ commute up to conjugacy on $M_c\cup\langle h,h'\rangle U$, 
contradicting the maximality of $M_c$. Note also that in this latter case the 
stabiliser of each component of $\langle h,h'\rangle U$ is reduced to the 
identity which clearly extends to $\langle h,h'\rangle U$.

\medskip

Assume we are in case (b). Consider the restriction of $h$ and
$h_\alpha=h{'}^{-\alpha}hh{'}^\alpha$ to $U$. Since $h$ and $h'$ commute on
$M_c$, $h$ and $h_\alpha$ coincide on $T$. Let $V$ be the geometric piece of 
the $JSJ$-decomposition for $M$ adjacent to $T$ and contained in $U$. Using 
Lemma \ref{lem:commutation}, we see that $h$ and $h_\alpha$ commute on $V$ and 
thus coincide on it, because they coincide on $T$. Thus $h$ and $h'$ commute on
$M_c\cup_{\alpha=0}^{p-1}h{'}^\alpha(V)$, and again we reach a contradiction to
the maximality of $M_c$.
\qed

\bigskip

We can thus assume to be in case (a) and that the $\langle h,h'\rangle$-orbit
of $T$ has $p$ elements.

\smallskip

\begin{Claim}\label{claim:winding number}
Each torus in the boundary of $M_c/<h>$ has winding number $p$ with respect
to $K$.
\end{Claim}

\demo
Since a boundary component $T$ of $M_c/<h>$ lifts to $p$ boundary components of 
$M_c$, the winding number of $T$ with respect to $K$ must be a multiple of $p$. 
We shall now reason by induction on  the number $n$ of boundary components of 
$M_c/<h>$. If $n=0$ there is nothing to prove.

If $n = 1$ the quotient spaces $M_c/<h>$ and $M_c/<h'>$ are solid tori, i.e. 
the exterior of a trivial knot which can be identified with a meridian of each 
solid torus. Note that the winding number of $T$ is precisely the linking 
number of $K$ with such a meridian. Note, moreover, that the spaces $M_c/<h>$ 
and $M_c/<h'>$ have a common quotient $\OO$ which is obtained by quotienting 
$M_c/<h>$ (respectively $M_c/<h'>$) via the the symmetry $\psi$ (respectively 
$\psi'$) of order $p$ and with non-empty fixed-point set, induced by $h'$ 
(respectively $h$). Since $\psi'$ preserves $\partial{(M_c/<h'>)}$ and has 
non-empty fixed-point set, $Fix(\psi')$ and the meridian of 
$\partial{(M_c/<h'>)}$ must form a Hopf link, in particular, their linking 
number is $1$. The image of $Fix(\psi')$ and of the meridian of 
$\partial{(M_c/<h'>)}$ form again a Hopf link in $\OO =(M_c/<h'>)/\psi$. By 
lifting them up to $M_c/<h>$ we see that the meridian lifts to a meridian and 
the image of $Fix(\psi')$ lifts to $K$ which thus have linking number $p$. 
Hence the property is proved in this case.

If $n>1$, we shall perform trivial Dehn surgery on $n-1$ boundary components of
$M_c/<h>$. Note that such a surgery does not change the winding number of the 
remaining boundary components (for the boundary components are unlinked), that 
the symmetry of order $p$ of $M_c/<h>$ extends to the resulting solid torus, 
and that the surgery can be lifted on $M_c$ in such a way that the quotient of 
the resulting manifold by the action of the diffeomorphism induced by $h'$ is 
again a solid torus. This last property follows from the fact that each 
connected component of $(E(K)\setminus(M_c/\langle h\rangle)$ is the exterior 
of a knot which lifts in $M$ to $p$ diffeomorphic copies. These $p$ copies of 
the knot exterior are permuted by $h'$ and a copy appears in the 
$JSJ$-decomposition of $E(K')$. This means that on each boundary component 
there is a well-defined meridian-longitude system which is preserved by $h$ and 
$h'$ and by passing to the quotient. The claim follows now from case $n=1$.
\qed

\bigskip

Now Claims \ref{claim:permutation} and \ref{claim:winding number} imply that 
$E_p(K)$ is a submanifold of $M_c/\langle h\rangle \cap E(K)$. 

Note, moreover, that the fixed-point set of the induced symmetry is contained 
in $M_f/\langle h\rangle\subset M_c/\langle h\rangle$. In particular, each 
torus of the $JSJ$-family separating such fixed-point set from $K$ lifts to a 
single torus of the $JSJ$-family for $M$ and its winding number cannot be a 
multiple of $p$. We can thus conclude that the fixed-point set of the symmetry 
induced by $h'$ is contained in $E_p(K)$. This finishes the proof of
Proposition \ref{prop:orbifold}.
\qed

\bigskip

\begin{Remark}\label{rem:orbifold} 
Note that $M_c/h \cap E(K)$ can be larger than $E_p(K)$ for there might be tori 
of the $JSJ$-collection for $M$ which have an $\langle h,h'\rangle$-orbit 
containing $p^2$ elements and which project to tori with winding number $p$. 
Note also that $E_p(K)$ coincides with $E(K)$ if there are no 
$JSJ$-tori in $E(K)$ with winding number $p$.
\end{Remark}

\medskip

\begin{Remark}\label{rem:commutation} 
The deck transformations $h$ and $h'$ cannot commute on the submanifolds $U$ of 
$M$ corresponding to branches of $\Gamma$ whose $h$- and $h'$-orbits coincide 
and consist of $p$ elements, if $h$ and $h'$ are different; that is, the 
stabiliser $h'h^{-1}$ is a finite order diffeomorphism of $U$ if and only if it 
is trivial. To see this, assume that there is a unique orbit of this type and 
assume by contradiction that $h$ and $h'$ commute on $M$ and are distinct. The 
diffeomorphism $h'$ would induce a non trivial symmetry of $E(K)$ of order $p$ 
and non-empty fixed-point set which fixes set-wise the projection of $U$ and 
acts freely on it. This contradicts the first part of Lemma 
\ref{lem:companion}. If there are $n>1$ such orbits an equivariant Dehn surgery 
argument on $n-1$ components leads again to a contradiction
\end{Remark}

\medskip

Here is a straightforward corollary of Proposition \ref{prop:orbifold} which 
generalises a result proved by B. Zimmermann \cite{Zim1} for hyperbolic knots.

\smallskip

\begin{Corollary}\label{cor:p-symmetry}
Let $K$ be a prime knot and let $p$ be an odd prime number. If $K$ has no
companion of winding number $p$ and has a $p$-twin, then $K$ admits a 
rotational symmetry of order $p$ with trivial quotient.
\qed
\end{Corollary}

\bigskip

So far we have proved that if a prime knot $K$ has a $p$-twin either $E(K)$ 
admits a $p$-rotational symmetry or a well-specified submanifold $E_p(K)$ of 
$E(K)$ admits a symmetry of order $p$ with non-empty fixed-point set. We shall 
say that the $p$-twin induces a \emph{symmetry}, respectively a \emph{partial 
symmetry}, of $K$.

\smallskip

\begin{Proposition}\label{prop:twin}
Let $K$ be a prime knot. Assume that $K$ has a $p$-twin and a $q$-twin for two 
distinct odd prime numbers.

\item{(i)} At least one twin, say the $q$-twin, induces a $q$-rotational  
symmetry $\psi_q$ of $K$. Moreover:

\item{(ii)} If the $p$-twin induces a partial $p$-symmetry of $K$, then 
$\partial E_p(K) \setminus \partial E(K)$ is a $JSJ$-torus which separates 
the fixed point set  $Fix(\psi_q)$ from $\partial E(K)$.
\end{Proposition}

\medskip

First we study some properties of partial symmetries induced by $p$-twins for 
an odd prime number $p$.

\smallskip

\begin{Lemma}\label{lem:partial companion}
Let $K$ be a prime knot and let $\psi$ be the partial symmetry of order $p$ 
induced on $E_p(K)$ by a $p$-twin. Let $T$ be a torus of the $JSJ$-collection 
of $E_p(K)$ which is not in the boundary. Then $T$ does not separate $\partial 
E(K)$ from $Fix(\psi)$ if and only if its $\psi$-orbit has $p$ elements. 
Moreover, this is the case if and only if the lift of $T$ to the $p$-fold 
cyclic branched cover of $K$ has $p$ elements.
\end{Lemma}

\demo
It suffices to perform $\psi$-equivariant Dehn fillings on the boundary
components $\partial E_p(K) \setminus \partial E(K)$ of $E_p(K)$ in such a way 
that the resulting manifold is a knot exterior $E(\hat K)$ and that the graph 
dual to the $JSJ$-decomposition of $E(\hat K)$ remains unchanged after filling 
(see the proof of Theorem \ref{thm:rotations}). Part (i) of Lemma 
\ref{lem:companion} then applies to the resulting knot $\hat K$ and the induced 
rotational symmetry. To apply Lemma \ref{lem:torus} it suffices to note that, 
as in the proof of Claim \ref{claim:winding number}, the fillings can be chosen 
in such a way that the induced fillings on the quotient $E_p(K)/\langle \psi 
\rangle$ give also a solid torus (see Remark \ref{rem:lift}).
\qed

\bigskip

\begin{Remark}\label{rem:case b} 
In particular, case (b) of the proof of Claim \ref{claim:permutation} cannot
happen for a torus $T$ in the situation of Lemma \ref{lem:partial companion}.
\end{Remark}

\medskip

\begin{Lemma}\label{lem:vertex}
Let $K$ be a prime knot and let $\psi$ be the partial symmetry of order $p$ 
induced on $E_p(K)$ by a $p$-twin. Let $T \subset \partial E_p(K) \setminus 
\partial E(K)$ be a torus which is $\psi$-invariant. Let $e_{T}$ be the 
corresponding edge in the tree dual to the $JSJ$-decomposition of $E_p(K)$. Let 
$v_K$ and $v_\psi$ be the vertices corresponding to the geometric pieces 
containing $\partial E(K)$ and $Fix(\psi)$ respectively. Then $v_\psi$ belongs 
to the unique geodesic joining $v_K$ to $e_{T}$ in this $JSJ$-tree.
\end{Lemma}

\demo 
If we cut $\S^3$ along a torus of the $JSJ$-collection of $E_p(K)$, the 
connected component which does not contain $K$ is a knot exterior and is thus 
contained in a ball in $\S^3$. If the conclusion of the Lemma were false, then 
we could find two tori of the $JSJ$-decomposition of $E_p(K)$ contained in two 
disjoint balls, one torus separating $Fix(\psi)$ from $K$ and the other 
coinciding with $T$ or separating it from $K$. In particular the linking number 
of $Fix(\psi)$ and a meridian of the solid torus bounded by $T$ (i.e. the 
winding number of $T$ with respect to $Fix(\psi)$) would be zero. This is 
impossible since $\psi$ leaves set-wise invariant $T$.
\qed

\bigskip

\begin{Remark}\label{rem:adjacency}
Lemma \ref{lem:vertex} has two interesting consequences. Since $h$ and $h'$
play symmetric roles, we deduce that $Fix(\psi)$ and $\partial E(K)$ must 
belong to the same geometric piece of the $JSJ$-decomposition of $E_p(K)$. This 
follows from the fact that, in $E_p(K') \cup \UU(K')$, $Fix(\psi)$ maps to 
$K'$, $K$ maps to $Fix(\psi')$, and $T$ maps to a $\psi'$-invariant torus. 
Moreover, each invariant boundary torus $T$ is adjacent to the geometric 
component containing $Fix(\psi)$ and $K$, else, we would get a contradiction to 
Lemma \ref{lem:partial companion}.
\end{Remark}

\medskip

{\bf Proof of Proposition \ref{prop:twin}(i).}
We argue by contradiction, assuming that there are a $p$-twin and a $q$-twin of 
$K$ which induce only partial symmetries of $E(K)$ for two distinct odd 
prime numbers $p$ and $q$. Then $\partial E_p(K)$ and $\partial E_q(K)$ are not 
empty. Moreover, we must have $E(K)\setminus E_p(K) \subset E_q(K)$ since the 
winding number along nested tori is multiplicative and thus the winding number 
of any $JSJ$-torus  contained in $E(K)\setminus E_p(K)$ must be of the form 
$kp$ and cannot be $q$. In particular $\partial E_p(K) \setminus \partial E(K) 
\subset \text{int}(E_q(K))$.

Let $T \in \partial E_p(K) \setminus \partial E(K)$ be a torus and let $\psi$ 
be the $q$-symmetry with non-empty fixed-point set induced on $E_q(K)$ by the 
$q$-twin. Since the winding number of $T$ is $p$, its lift to the $q$-fold 
cyclic branched cover of $K$ is connected. According to part (i) of Lemma 
\ref{lem:companion} and to Lemmata \ref{lem:torus} and \ref{lem:partial 
companion}, $T$ must separate $\partial E(K)$ from $Fix(\psi)$. Since 
$Fix(\psi)$ is connected, we see that so must be $\partial E_p(K) \setminus 
\partial E(K) = T$. The final contradiction is then reached by applying Remark 
\ref{rem:adjacency}.
\qed 

\bigskip

{\bf Proof of Proposition \ref{prop:twin}(ii).} 
This is a consequence of the proof of part (i): note that in the proof $\psi$ 
may be a global or partial symmetry.
\qed

\bigskip

We are now in a position to prove Theorem \ref{thm:twins}.

\medskip

{\bf Proof of part (i) of Theorem \ref{thm:twins}.}
We argue by contradiction, assuming that $K$ admits twins for three distinct, 
odd prime numbers $p, q, r$. Under this assumption, it follows that $K$ is a 
non-trivial knot. 

If the three twins induce rotational symmetries of the knot $K$, then part (i) 
of Theorem \ref{thm:rotations} gives a contradiction.

Therefore part (i) of Proposition \ref{prop:twin} implies that twins of orders, 
say $q$ and $r$, induce rotational symmetries $\psi_q$ and $\psi_r$ of $K$ 
having order $q$ and $r$ respectively, while a $p$-twin induces only a partial 
rotational symmetry of $E(K)$ of order $p$. 

Then part (ii) of Proposition \ref{prop:twin} shows that $\partial E_p(K) 
\setminus \partial E(K)$ is a $JSJ$-torus in $E(K)$ which separates $\partial 
E(K)$ from both $Fix(\psi_q)$ and $Fix(\psi_r)$. This contradicts part (ii) of 
Theorem \ref{thm:rotations} which states that $Fix(\psi_q)$ and 
$Fix(\psi_r)$ must sit in the $JSJ$-component containing $\partial E(K)$.
\qed

\bigskip

{\bf Proof of part (ii) of Theorem \ref{thm:twins}.}
Let $K$ be a prime knot and let $p$ be an odd prime number. We assume that $K$ 
has at least two non-equivalent $p$-twins $K_1$ and $K_2$ and look for a 
contradiction. 

\medskip

If both  $\psi_{1}$ and $\psi_{2}$ are rotational symmetries of order $p$ of
$K$, then by M. Sakuma \cite[Thm. 3]{Sa2} they are conjugate since $K$ is 
prime. This would contradict the hypothesis that the knots $K_1$ and $K_2$ 
are not equivalent.

\medskip

Assume now that at least one symmetry, say $\psi_{1}$ is partial. Then 
$\psi_{1}$ and $\psi_{2}$ are rotational symmetries of order $p$ of the 
submanifold $E_p(K) \subset E(K)$. Let $X_0$ be the geometric piece of 
the $JSJ$-decomposition of $E(K)$ containing $\partial E(K)$. Then $\psi_1$ 
(respectively $\psi_2$) generates a finite cyclic subgroup $G_1$ (respectively
$G_2$) of the group $Diff^{+,+}(X_0, \partial E(K))$ of diffeomorphisms of the 
pair $(X_0, \partial E(K))$ which preserve the orientations of $X_0$ and of 
$\partial E(K)$. Moreover, one can assume that $G_1$ and $G_2$ act
geometrically on $X_0$.

If $X_0$ admits a hyperbolic structure, it is a consequence of the proof of the 
Smith conjecture (see for example \cite[Lemma 2.2]{Sa2}) that the subgroup of 
$Diff^{+,+}(X_0, \partial E(K))$ consisting of restrictions of isometries of 
$X_0$ is finite cyclic. Hence $G_1 = G_2$ and up to taking a power 
$\psi_1 = \psi_2$ on $X_0$.

If $X_0$ is Seifert fibred, then it must be a cable space, since $K$ is prime. 
The uniqueness of the Seifert fibration and the fact that the basis of the 
Seifert fibration has no symmetry of finite order imply that the cyclic groups 
$G_1$ and $G_2$ belong to the circle action $S^1 \subset Diff^{+,+}(X_0, 
\partial E(K))$ inducing the Seifert fibration of $X_0$, see 
\cite[Lemma 2.3]{Sa2}. Since $G_1$ and $G_2$ have the same prime order, up to 
taking a power $\psi_1 = \psi_2$ on $X_0$.

Let $h_1$ and $h_2$ be the deck transformations on $M$ associated to the
$p$-fold cyclic coverings branched along $K_1$ and $K_2$, and which induce 
$\psi_1$ and $\psi_2$. Then by taking a suitable powers, $h_1$ and $h_2$ 
coincide up to conjugacy on the geometric piece $\widetilde X_0$ of the 
$JSJ$-decomposition of $M$ containing the preimage of $K$. The following lemma 
shows that they will coincide on $M$, contradicting our hypothesis.
\qed

\bigskip

\begin{Lemma}\label{cor:extension}
If the covering transformations $h$ and $h'$ preserve a $JSJ$-piece or a
$JSJ$-torus of $M$ and coincide on it, then they can be chosen, up to
conjugacy, to coincide everywhere.
\end{Lemma}

\demo 
This is  a consequence of the proofs of Propositions \ref{prop:subtree} and
\ref{prop:orbifold}. We shall start by showing that we can always assume that 
there is a piece $V$ of the $JSJ$-decomposition on which $h$ and $h'$ coincide. 
To this purpose, assume that $h$ and $h'$ coincide only on a $JSJ$-torus $T$. 
According to Lemma \ref{lem:commutation} and Remark \ref{rem:same action}, $h$ 
and $h'$ coincide on the geometric pieces of the decomposition adjacent to $T$, 
which are also invariant. Consider now the maximal subtree $\Gamma_1$ of 
$\Gamma$ such that the restrictions of $h$ and $h'$ to the corresponding 
submanifold $M_1$ of $M$ coincide, up to conjugacy, and such that 
$V\subset M_1$. Let $S$ be a $JSJ$-torus for $M$ in the boundary of $M_1$. 
Since $h$ and $h'$ coincide on $M_1$, the $h$-orbit and the $h'$-orbit of $S$ 
coincide as well and consist of either one single element $\{S\}$ or $p$ 
elements $\{S,h(S)=h'(S),...,h^{p-1}(S)={h'}^{p-1}(S)\}$. In the former case, 
according to Lemma \ref{lem:commutation}, $\Gamma_1$ would not be maximal. In 
the latter case, we are precisely in the situation described in part (a) of 
Claim \ref{claim:permutation}. Once more, $\Gamma_1$ is not maximal because one 
can impose that $h$ and $h'$ act in the same way on the $p$ connected 
components with connected boundary obtained by cutting $M$ along the $\langle 
h,h'\rangle$-orbit of $S$ (see Remark \ref{rem:commutation}). This 
contradiction shows that $M=M_1$ and the lemma is proved.
\qed

\bigskip

{\bf Proof of part (iii) of Theorem \ref{thm:twins}.}
First we analyse the case of a knot admitting two twins, one of which induces a 
partial symmetry. 

\smallskip

\begin{Proposition}\label{prop:partial}
Let $K$ be a prime knot admitting a $p$-twin $K'$ and a $q$-twin $K''$ for two distinct
odd prime numbers $p$ and $q$. If $K'$ induces a partial symmetry of $K$ then
$K'$ and $K''$ are not equivalent.
\end{Proposition}

\demo
By part (ii) of Proposition \ref{prop:twin}, $E_p(K)$ has a unique boundary 
component which separates $\partial E(K)$ from the fixed-point set of the
$q$-rotational symmetry $\psi$ induced by $K''$. By cutting $\S^3$ along 
$T = \partial E_p(K)$ we obtain a solid torus $V=E_p(K)\cup \UU(K)$ containing 
$K$, and a knot exterior $E_T$. $K$ admits a $q$-rotational symmetry $\psi$ 
induced by $K''$ which preserves this decomposition and induces a 
$q$-rotational symmetry with trivial quotient (see Lemma \ref{lem:companion}) 
on $E_T$ and a free $q$-symmetry $\tilde\psi$ on $V$. The covering 
transformation for the knot $K'$ induces a $p$-symmetry $\varphi$ of $V$ with 
non-empty fixed-point set.

Assume now by contradiction that $K'=K''$. Since $K'$ induces a partial 
symmetry of $K$ and vice versa, $S^3$ admits a decomposition into two 
pieces: $V'=E_p(K')\cup \UU(K')$ and $E_T $. On the other hand, since $K''$ 
induces a genuine $q$-rotational symmetry of $K$, $K''$ admits a 
$q$-rotational symmetry $\psi''$ induced by $K$ which preserves the 
aforementioned decomposition and induces a $q$-rotational symmetry with trivial 
quotient on $E_T$. Using the fact that $E_T$ is the exterior of a prime knot 
(see Lemma \ref{lem:prime}) and M. Sakuma's result \cite[Thm. 3]{Sa2}, we see 
that the two $q$-rotational symmetries with trivial quotient induced by $\psi$ 
and $\psi''$ on $E_T$ act in the same way. Let now $E_0$ be the smallest knot 
exterior of the $JSJ$-decomposition of $E_T$ on which $\psi=\psi''$ induces a 
$q$-rotational symmetry with trivial quotient (this is obtained by cutting 
$E_T$ along the torus of the $JSJ$-decomposition closest to $Fix(\psi)$ 
-respectively $Fix(\psi'')$- and separating it from $T$. Consider now the lift,
denoted by $(X,\KK)$, to $(S^3,K'')$ of $(E_0, Fix(\psi))/\psi$. We claim that 
$(X,\KK)=(V',K')$. Indeed, $X$ contains $K''=K'$ by construction, and its 
boundary is the unique torus of the $JSJ$-decomposition which is left invariant 
by the $q$-rotational symmetry of $K''$ -by construction again- and which is 
closest to $K''$ (compare Remark \ref{rem:adjacency}). Since 
$E_0/\psi=E_0/\psi''$, and a solid torus has a unique $q$-fold cyclic cover, we 
deduce that $(V',K')=(X,\KK)=(V,K)$. In particular, the deck transformations 
for $K$ and $K'$ on their common $p$-fold cyclic branched cover can be chosen 
to coincide on the lift of $V=V'$. Lemma \ref{cor:extension} implies that 
$K=K'$ contradicting the fact that $K'$ is a $p$-twin.
\qed

\bigskip

Let $K'$ be a $p$-twin and a $q$-twin of $K$ for two distinct odd prime numbers 
$p$ and $q$. Proposition \ref{prop:partial} implies that $K'$ induces two 
rotational symmetries $\psi_p$ and $\psi_q$ of $K$ with trivial quotients and 
orders $p$ and $q$. Part (ii) of Theorem \ref{thm:rotations} shows that the 
fixed-point sets $Fix(\psi_p)$ and $Fix(\psi_q)$ lie in the $JSJ$-component of 
$E(K)$ which contains $\partial E(K)$. Then the proof of part (iii) of Theorem 
\ref{thm:twins} follows from the following:

\smallskip

\begin{Lemma}\label{lem:commuting symmetries}
Let $K$ be a prime knot admitting two rotational symmetries $\psi$ and
$\varphi$ of odd prime orders $p > q$. If the fixed-point sets of $\psi$ and 
$\varphi$ lie in the component which contains $\partial E(K)$, then the two 
symmetries commute up to conjugacy.
\end{Lemma}

\demo
Reasoning as in the proof of part (ii) of Theorem \ref{thm:twins}, one can show
that $\psi$ and $\varphi$ commute on the component which contains $\partial
E(K)$. Since all other components are freely permuted according to part (i) of
Lemma \ref{lem:companion}, the conclusion follows as in the proof of part (a)
of Claim \ref{claim:permutation}.
\qed

\bigskip

{\bf Proof of Corollary \ref{cor:composite}.}
First of all note that, because of the uniqueness of the Milnor-Kneser 
decomposition of the covers of $K$ and $K'$, the number of prime summands of $K$ 
and $K'$ is the same. After ditching components of $K$ and $K'$ that appear in 
both decompositions in equal number, we can assume that $K_i$ is not equivalent 
to $K'_\ell$, for all $i,\ell=1,...,t$. If $K$ and $K'$ have three common 
cyclic branched covers of odd prime orders, we deduce that for each 
$i=1,...,t$, $K_i$ is not determined by its $p_j$-fold cyclic branched cover, 
$j=1,2,3$, for it is also the $p_j$-fold cyclic branched cover of some 
$K'_{i_j}$ not equivalent to $K_i$. Hence $K_i$ would have twins for three
distinct odd prime orders which is impossible by Theorem \ref{thm:twins}. 
\qed

\section{Examples}

Examples of prime knots admitting a $p$-twin which induces a global rotational
symmetry of order $p$ were first constructed by Y. Nakanishi \cite{Na} and M.
Sakuma \cite{Sa1}. They considered a prime link with two trivial components 
whose linking number is $1$. By taking the $p$-fold cyclic cover of $\S^3$ 
branched along the first (respectively the second) component of the link one
gets again $\S^3$ and the second (respectively first) component lifts to a
prime knot. The two knots thus constructed have the same $p$-fold cyclic
branched cover by construction (see also Remark \ref{rem:lift}), moreover, by 
computing their Alexander polynomial they were shown to be distinct.

In \cite[Thm 3 and Cor. 1]{Zim1} B. Zimmerman showed that if a hyperbolic knot 
has a $p$-twin, for $p\ge 3$, then the $p$-twin induces a global symmetry and 
the two knots are thus obtained by Y. Nakanishi and M. Sakuma's construction 
where the quotient link is hyperbolic and admits no symmetry which exchanges 
its two components.

As a matter of fact, the links considered by Y. Nakanishi and M. Sakuma are in 
fact hyperbolic and so are the resulting twins if $p$ is at least $3$, 
according to the orbifold theorem \cite{BoP}, see also \cite{CHK}. Note that, 
when $p=2$, the situation, even in the case of hyperbolic knots, is much more 
complex and there are several ways to construct $2$-twins of a given knot. In 
this section we shall see how one can construct, for each given odd prime $p$, 
two prime, non simple, knots which are $p$-twins, and such that the symmetries 
they induce are not global.

The first construction shows that the number $\nu$ of components of $\partial 
E_p(K)\setminus\partial E(K)$ can be arbitrarily large. This means that the 
situation encountered in Proposition \ref{prop:twin}(ii) is extremely special. 
The second construction shows that our result is indeed best possible even for
prime knots with $p$-twins inducing partial symmetries: we shall construct 
prime knots admitting a $p$-twin inducing a partial symmetry and a $q$-twin 
inducing a global rotational symmetry. 

\bigskip

\subsection{Knots admitting a $p$-twin inducing only a partial symmetry} 

\medskip

Assume we are given a hyperbolic link $L=L_1\cup...\cup L_{\nu+2}$, with $\nu+2
\ge 3$ components, satisfying the following requirements:

\medskip

{\bf Property $*$}

\begin{enumerate}

\item The sublink $L_3\cup...\cup L_{\nu+2}$ is the trivial link;

\item For each $i=1,2$ and $j=3,...,\nu+2$, the sublink $L_i\cup L_j$ is a Hopf
link;

\item ${\rm lk}(L_1,L_2)$ is prime with $p$;

\item No symmetry of $L$ exchanges $L_1$ and $L_2$.

\end{enumerate}

\medskip
We shall consider the orbifold $\OO=(\S^3,(L_1\cup L_2)_p)\setminus
\UU(L_3\cup...\cup L_{\nu+2})$ which is the $3$-sphere with singular set of
order $p$ the (sub)link $L_1\cup L_2$ and an open tubular neighbourhood of the
(sub)link $L_3\cup...\cup L_{\nu+2}$ removed. $\OO$ is hyperbolic if $p\ge3$, 
and will represent the quotient of $\OO_p=E_p(K)\cup\UU(K)$ and 
$\OO_p'=E_p(K')\cup\UU(K')$ via the action of the partial $p$-symmetries. 
Indeed, to obtain $\OO_p$ (respectively $\OO_p'$) take the $p$-fold cyclic 
orbifold cover of $(\S^3,(L_1\cup L_2)_p)\setminus
\UU(L_3\cup...\cup L_{\nu+2})$ which desingularises $L_2$ (respectively $L_1$).
Observe that one can fix a longitude-meridian system on each boundary 
component of $\OO$, induced by those of $L_i$, $i=3,\dots,\nu+2$. Note that, 
because of condition 4 of Property $*$, the two orbifolds $\OO_p$ and $\OO_p'$ 
with the fixed peripheral systems are distinct. 

Remark that $\OO_p$ and $\OO'_p$ can be obtained by the orbifold covers,
analogous to those described above, of $(\S^3,(L_1\cup L_2)_p)$ (which are 
topologically $\S^3$) by removing open regular neighbourhoods of the lifts of 
the components $L_3\cup...\cup L_{\nu+2}$. Note that these components 
lift to trivial components whose linking number with the lift of $L_i$, 
$i=1,2$, is precisely $p$, because of condition 2, and which form again a 
trivial link.

For each $j=3,...,\nu+2$, choose a knot exterior $E(\KK_j)$ to be glued along 
the $j$-th boundary component of $\OO_p$ and $\OO'_p$ in such a way that a 
fixed longitude-meridian system on $E(\KK_j)$ is identified with the lift of 
the longitude-meridian system on the $j$-th boundary component of $\OO$. The underlying spaces of the
orbifolds $\OO_p\cup_{j=3}^{\nu+2} E(\KK_j)$ and $\OO'_p\cup_{j=3}^{\nu+2} 
E(\KK_j)$ are topologically $\S^3$ and it is easy to see that their singular 
sets are connected (see condition 3). The resulting knots have the same 
$p$-fold cyclic branched cover, however, since $\OO_p$ and $\OO'_p$ are 
distinct, they are not equivalent.

\bigskip

\begin{Remark}
Observe that we have just shown that the number of connected components of
$\partial E_p(K)\setminus\partial E(K)$, which is precisely $\nu$, can be 
arbitrarily large. Note also that if $\nu\ge2$, according to Proposition
\ref{prop:twin}, the knot $K$ has no $q$-twins for $q\neq p$ odd prime.
\end{Remark}

\bigskip

We shall now prove that links with Property $*$ exist. Notice that for $\nu=1$
links satisfying all the requirements where constructed by Zimmermann in
\cite{Zim2}, see also \cite{Pao1}.

\begin{figure}
\begin{center}
\includegraphics[width=8cm]{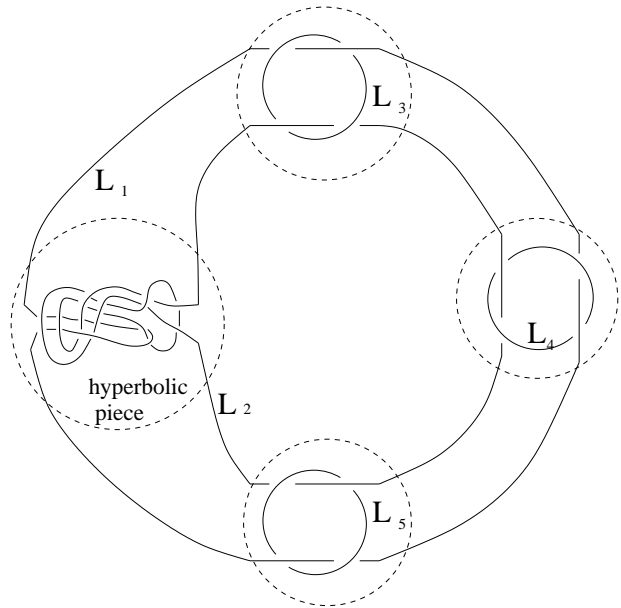}
\end{center}
\caption{The link $L$ and its Bonahon-Siebenmann decomposition.}
\label{fig:one}
\end{figure}

Consider the link given in Figure \ref{fig:one} for $\nu=3$ (the generalization 
for arbitrary $\nu\ge1$ is obvious). Most conditions are readily checked just 
by looking at the figure, and we only need to show that $L$ is hyperbolic and 
has no symmetries which exchange $L_1$ and $L_2$. To this purpose, we shall 
describe the Bonahon-Siebenmann decomposition of the orbifold $(\S^3,(L)_2)$, 
where all components have $\Z/2\Z$ as local group. The decomposition consists 
of one single hyperbolic piece (see Figure \ref{fig:one}) and $\nu+1$ 
(respectively $1$) Seifert fibred pieces if $\nu\ge2$ (respectively $\nu=1$). 
Since the Seifert fibred pieces contain no incompressible torus, the 
hyperbolicity of $L$ follows.

Note now that every symmetry of $L$ must leave invariant the unique hyperbolic
piece of the decomposition. This piece is obtained by quotienting the
hyperbolic knot $10_{155}$ via its full symmetry group $\Z/2\Z\oplus\Z/2\Z$ and 
thus has no symmetries (for more details see \cite{Pao1}), so we conclude that 
the components $L_1$ and $L_2$ are non exchangeable.

\bigskip

\subsection{Knots admitting a $p$-twin inducing a partial symmetry and a
$q$-twin inducing a global symmetry}

\medskip

Let $\KK$ be a hyperbolic knot admitting a $p$-twin and a $q$-twin; the twins 
of $\KK$ induce global symmetries, so that $\KK$ admits a $p$- and a
$q$-rotational symmetry with trivial quotient (see \cite{Zim2}, where a method
to construct hyperbolic knots with two twins is described). Remove a tubular
neighbourhood of the axis of the symmetry of order $q$ (note that the two
symmetries have disjoint axes), and use the resulting solid torus $V$ to 
perform Dehn surgery on the exterior $E$ of the $(2,q)$-torus knot. Denote by 
$K$ the image of $\KK$ after surgery. We require that:

\begin{enumerate}

\item The resulting manifold is $\S^3$;

\item The $q$-rotational symmetry of $E$ and the restriction of the
$q$-rotational symmetry of $\KK$ to $V$ give a global $q$-rotational symmetry 
of $K$;

\item The $q$-rotational symmetry of $K$ has trivial quotient.

\end{enumerate}

\begin{figure}
\begin{center}
\includegraphics[width=10cm]{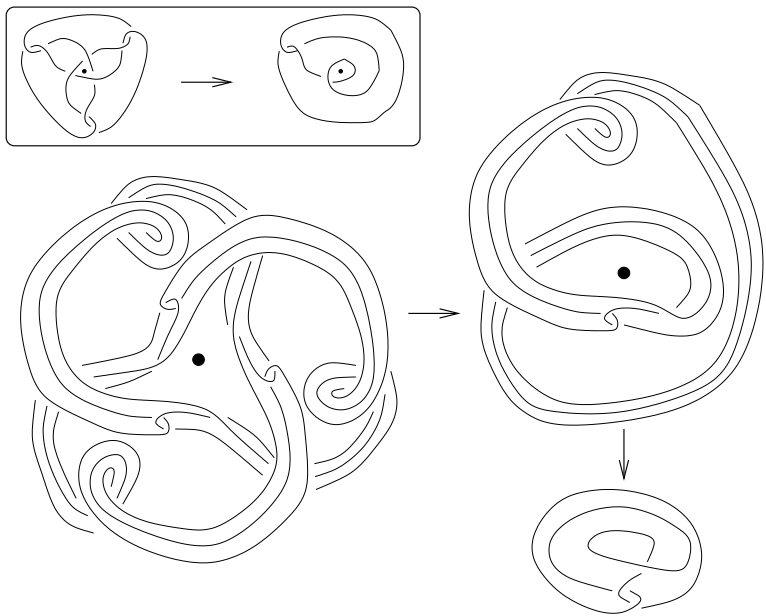}
\end{center}
\caption{Satellising so that the induced rotation has trivial quotient.}
\label{fig:two}
\end{figure}

Note that the last requirement can be met by choosing appropriately the
longitude when satellising, as illustrated in Figure \ref{fig:two}. We claim 
that $K$ admits a $q$-twin, $K''$, and a $p$-twin, $K'$. $K''$ is obtained by 
the standard method described in Remark \ref{rem:lift}. Note that $K\neq K''$, 
for the roots of the $JSJ$-decompositions of the exteriors of $K$ and $K''$ are 
hyperbolic and Seifert fibred respectively. To construct $K'$, consider the 
$p$-twin $\KK'$ of $\KK$ and let $V'$ be the solid torus obtained by removing 
the axis of the $q$-rotational symmetry of $\KK'$. Note that $V$ and $V'$ have 
a common quotient obtained by taking the space of orbits of the $p$-rotational 
symmetries, however $V$ and $V'$ are different orbifolds by construction. Fix a 
longitude-meridian system on $V$ (the one used for the surgery): by first 
quotienting and then lifting it, get a longitude-meridian system on $V'$ that 
must be used to perform surgery along a copy of $E$. The image of $\KK'$ after 
the surgery will be $K'$. Note that, when taking the $p$-fold cyclic branched 
covers of $K$ and $K'$, the hyperbolic orbifolds $V$ and $V'$ lift to the same 
manifold by construction, while the Seifert fibred part lifts, in both cases, 
to $p$ copies of $E$. Again by construction, the gluings are compatible and the 
two covers coincide. It is also evident that $K'$ can only induce a partial 
symmetry of $K$, and the claim is proved.

\bigskip

\begin{Remark}
Note that according to Proposition \ref{prop:partial} the $p$-twins and 
$q$-twins obtained in this construction cannot be equivalent.
\end{Remark}

\section{Homology spheres as cyclic branched covers}

By the proof of the Smith conjecture Corollary \ref{cor:homologysphere} is true 
for the $3$-sphere $S^3$. So from now on we assume that the integral homology 
sphere $M$ is not homeomorphic to $S^3$. Then by \cite[Thm1]{BPZ}, $M$ can be a 
$p_i$-fold cyclic branched cover of $\S^3$ for at most three pairwise distinct 
odd prime numbers $p_i$. Moreover if $M$ is irreducible and is the $p_i$-fold 
cyclic branched cover of $\S^3$ for three pairwise distinct odd prime numbers 
$p_i$, then the proof of \cite[Corollary 1.(i)]{BPZ} shows that for each prime 
$p_i$, $M$ is the $p_i$-fold cyclic branched cover of precisely one knot. Since 
a knot admits at most one $p$-twin for an odd prime integer $p$, we need only 
to consider the case when the irreducible integral homology sphere $M$ is the 
branched cover of $\S^3$ for precisely two distinct odd primes, say $p$ and 
$q$. Moreover \cite[Corollary 1.(ii)]{BPZ} shows that $M$ has a non trivial 
$JSJ$-decomposition. 

Looking for a contradiction, we can assume that, for each prime, $M$ is the 
branched covering of two distinct knots with covering transformations $\psi$, 
$\psi'$ of order $p$ and $\varphi$, $\varphi'$ of order $q$. 

If each rotation of order $p$ commutes with each rotation of order $q$ up to 
conjugacy, then the contradiction follows from the following claim which is an 
easy consequence of Sakuma's result \cite[Thm. 3]{Sa2} (see \cite[Claim 
8]{BPZ}). 

\begin{Claim}\label{claim:unique symmetry} 
Let $n\ge3$ be a fixed odd integer. Let $\rho$ be a rotation with trivial 
quotient of an irreducible manifold $M$. All the rotations of $M$ of order $n$ 
which commute with $\rho$ are conjugate in $Diff(M)$ into the same cyclic group 
of order $n$.  
\qed
\end{Claim}

Otherwise, consider the subgroup $G=\langle \psi, \psi', \varphi, \varphi' 
\rangle$ of diffeomorphisms of $M$. According to the proof of \cite[Proposition 
4]{BPZ}, each rotation of order $p$ commutes with each rotation of order $q$ up 
to conjugacy, unless the induced action of $G$ on the dual tree of the 
$JSJ$-decomposition for $M$ fixes precisely one vertex corresponding to a 
hyperbolic piece $V$ of the decomposition and $\{p,q\}=\{3,5\}$. In this case, 
one deduces as in the proof of \cite[Corollary 1.(ii)]{BPZ} that the 
restrictions of $\psi$ and $\psi'$ (respectively $\varphi$ and $\varphi'$) 
coincide up to conjugacy on $V$. Then the desired contradiction follows from 
Lemma \ref{cor:extension} which implies that $\psi$ and $\psi'$ (respectively 
$\varphi$ and $\varphi'$) coincide up to conjugacy on $M$.
\qed

\begin{footnotesize}

\bigskip

{\textsc{Laboratoire Emile Picard (UMR 5580 du CNRS)}}

{\textsc{Universit\'e Paul Sabatier}}

{\textsc{118 route de Narbonne}}

{\textsc{31062 Toulouse cedex 4, France}}

{\tt{boileau@picard.ups-tlse.fr}}

\medskip

{\textsc{I.M.B.(UMR 5584 du CNRS)}}

{\textsc{Universit\'e de Bourgogne}}

{\textsc{B.P. 47 870}}

{\textsc{9 av. Alain Savary}}

{\textsc{21078 Dijon cedex, France}}

{\tt{paoluzzi@u-bourgogne.fr}}

\end{footnotesize}

\end{document}